\input amstex
\documentstyle{amsppt}
\document
\topmatter
\title
A rigid analytical regulator for the $K_2$ of Mumford curves
\endtitle
\author
Ambrus P\'al
\endauthor
\date
April 25, 2009.
\enddate
\address Department of Mathematics, 180 Queen's Gate, Imperial College,
London SW7 2AZ, United Kingdom\endaddress
\email a.pal\@imperial.ac.uk\endemail
\abstract We construct a rigid analytical regulator for the $K_2$ of
Mumford curves, a non-archimedean analogue of the complex analytical Beilinson-Bloch-Deligne regula\-tor.
\endabstract
\footnote" "{\it 2000 Mathematics Subject Classification. \rm 19C20
(primary), 11G20 (secondary).}
\endtopmatter

\heading 1. Introduction
\endheading

\definition{Motivation 1.1} Suitably normalized, the Beilinson-Bloch-Deligne regulator is a morphism of functors (for definition and properties see [7], pages 18-23):
$$\{\cdot\}:H^2_{\Cal M}(U,\Bbb Z(2))\rightarrow H_{an}^1(U,\Bbb C^*)=
\text{Hom}(\pi^{ab}_1(U),\Bbb C^*),$$
defined on the category of Riemann surfaces such that
\roster
\item"$(i)$" if $U=Y-\{x\}$, where $Y$ is a Riemann surface, $x\in Y$ is a point, and $f$, $g\in H^0(U,\Cal O^*)$ are meromorphic at $x$, then
the value of $\{f\otimes g\}$ on the positively oriented loop around $x$ is the tame symbol $\{f,g\}_x$ at $x$,
\item"$(ii)$" If $f_n\rightarrow f$ and $g_n\rightarrow g$ converge uniformly on compact sets, where $f_n$, $f$, $g_n$ and $g$ are elements of
$H^0(U,\Cal O^*)$, then the value of $\{f_n\otimes g_n\}$ on every closed
loop in $U$ converges to the value of $\{f\otimes g\}$ on that loop.
\endroster
The second property follows from the integral representation of the
monodromy. We wish to construct a regulator which is analogous to the
the monodromy of the Beilinson-Bloch-Deligne regulator on the boundary
components of Riemann surfaces with boundary in the rigid analytical
context. If we want the two properties above to hold, then we should
define the regulator of two nowhere vanishing holomorphic functions $f$
and $g$ on a connected rational subdomain of the projective line around
a complementary disk by approximating $f$ and $g$ by rational functions
and set the regulator $\{f,g\}$ as the limit of the product of tame
symbols at points inside the disk. This is exactly what we will do. As an application we will construct a rigid analytical regulator for the $K_2$ of
Mumford curves whose properties strongly resemble those of the complex analytical Beilinson-Bloch-Deligne regulator. It is a homomorphism:
$$\{\cdot\}:H^2_{\Cal M}(X,\Bbb Z(2))\rightarrow\Cal H(\Gamma(\goth X_0),\Bbb C^*),$$
where $\Bbb C$ is an algebraically closed field complete with respect to an ultrametric
absolute value, $X$ is a Mumford curve over $\Bbb C$ and $\Cal H(\Gamma(\goth X_0),\Bbb C^*)$
denotes the group of $\Bbb C^*$-valued harmonic cochains on the oriented incidence graph of the
special fiber $\goth X_0$ of some semi-stable model $\goth X$ of $X$ over the valuation ring
of $\Bbb C$. The problem of finding rigid analytic analogues of the Beilinson-Bloch-Deligne
regulator has already been studied for example in [3] by Coleman. Coleman's methods are
radically different from ours, and they cover different ground. On the other hand in the two
papers [8] and [9] Kato constructed a regulator for higher local fields which is essentially
the same as ours. Partially to generalize his results to fields whose valuation is not
discrete, partially to be self-contained, we develop the foundations of this theory
independently of his work.
\enddefinition
\definition{Contents 1.2} In the second chapter we define the regulator first for connected
rational subdomains of the projective line by approximating holomorphic functions on the domain
 by rational functions and taking the limit of the product of the tame symbols inside a
complementary open disk as we already mentioned above. The main result of this section is
Theorem 2.2 which gives a complete characterization of our regulator. In the third chapter we
first extend the definition of the rigid analytical regulator for the $K_2$ of connected
rational subdomains of the projective line (Theorem 3.2). Then we prove the invariance theorem
(Theorem 3.11) which is the functoriality property of the regulator for such subdomains. In
order to do so we develop an elementary homology theory for connected rational subdomains of
the projective line (Theorems 3.6 and 3.8). The fourth chapter is somewhat technical: it
contains two auxiliary results used in the next chapter. In the fifth chapter we reap the
fruits of our labours when we define the rigid analytical regulator for Mumford curves. The
latter takes values in the group of $\Bbb R^*$-valued harmonic cochains on the oriented
incidence graph of the special fiber of a semi-stable model of the Mumford curve where
$\Bbb R$ is the field of definition. We also show a reciprocity law relating the valuation of
the rigid analytical regulator to the (generalized) tame symbol along the special fiber. We
also formulate a functoriality property for this regulator phrased in terms of measures on the
ends of the universal covering of the oriented incidence graph of the special fiber
(Proposition 5.6). In the last chapter we look at the particular case of the Drinfeld upper
half plane and we express the tame symbol at the cusps as a non-archimedean integral of the
rigid analytic regulator (Theorem 6.5).
\enddefinition
\definition{Notation 1.3} In this paper we will use the somewhat incorrect notation $K_2(X)$
to denote $H^2_{\Cal M}(X,\Bbb Z(2))$ for various types of spaces $X$ as the latter is rather
awkward.
\enddefinition
\definition{Acknowledgement 1.4} I wish to thank the CRM and the IH\'ES, where this article was
written, for their warm hospitality and the pleasant environment they created for productive
research.
\enddefinition

\heading 2. The rigid analytical regulator for connected rational
subdomains
\endheading

\definition{Notation 2.1} Let $\Bbb C$ be an algebraically closed field complete with respect to an ultrametric absolute value which will be denoted by $|\cdot|$. Let $|\Bbb C|$ denote the set of values of the latter. Let $\Bbb P^1$  denote the projective line over $\Bbb C$. For any $x\in\Bbb P^1$ and any two rational non-zero functions $f$, $g\in\Bbb C((t))$ on the projective line let $\{f,g\}_x$ denote the tame symbol of the pair $(f,g)$ at $x$. We call a set $D\subset\Bbb P^1$ an open disk if it is the image of the set $\{z\in\Bbb C||z|<1\}$ under a M\"obius transformation. Recall that a subset $U$ of $\Bbb P^1$ is a connected rational subdomain, if it is non-empty and it is the complement of the union of finitely many pair-wise disjoint open disks. Let $\partial U$ denote the set of these complementary open disks. The elements of $\partial U$ are called the boundary components of $U$, by slight abuse of language. Let $\Cal O(U)$, $\Cal R(U)$ denote the algebra of holomorphic functions on $U$ and the subalgebra of restrictions of rational functions, respectively.  Let $\Cal O^*(U)$, $\Cal R^*(U)$ denote the group of invertible elements of these algebras. The group $\Cal R^*(U)$ consists of rational functions which do not have poles or zeros lying in $U$. For each $f\in\Cal O(U)$ let $\|f\|$ denote $\sup_{z\in U}|f(z)|$. This value is finite, and makes $\Cal O(U)$ a Banach algebra over $\Bbb C$. We say that the sequence $f_n\in\Cal O(U)$ converges to $f\in\Cal O(U)$, denoted by $f_n\rightarrow f$, if $f_n$ converges to $f$ with respect to  the topology of this Banach algebra, i.e.
$\lim_{n\rightarrow\infty}\|f-f_n\|=0$. For every real number
$0<\epsilon<1$ we define the sets $\Cal O_{\epsilon}(U)=\{f\in\Cal
O(U)|\|1-f\|\leq\epsilon\}$, and $U_{\epsilon}=\{z\in\Bbb C|
|1-z|\leq\epsilon\}$.
\enddefinition
The main result of this section is the following
\proclaim{Theorem 2.2} There is a unique map $\{\cdot,\cdot\}_D:\Cal
O^*(U)\times\Cal O^*(U)\rightarrow\Bbb C^*$ for every $D\in\partial U$,
called the rigid analytic regulator, with the following properties:

\noindent (i) For any two $f$, $g\in\Cal R^*(U)$ their regulator is:
$$\{f,g\}_D=\prod_{x\in D}\{f,g\}_x,$$

\noindent (ii) the regulator $\{\cdot,\cdot\}_D$ is bilinear
in both variables,

\noindent (iii) the regulator $\{\cdot,\cdot\}_D$ is alternating:
$\{f,g\}_D\cdot\{g,f\}_D=1$,

\noindent (iv) if $f$, $1-f\in\Cal O(U)^*$, then $\{f,1-f\}_D$ is $1$,

\noindent (v) for each $f\in\Cal O_{\epsilon}(U)$ and $g\in\Cal O^*(U)$
we have $\{f,g\}_D\in U_{\epsilon}$.
\endproclaim
\definition{Remark 2.3} It is an immediate consequence of property $(v)$
that the rigid analytic regulator is continuous with respect to  the
supremum topology. Explicitly, if $f$ and $g$ are elements of $\Cal
O^*(U)$, $D\in\partial U$ is a boundary component, and $f_n\in\Cal
O^*(U)$, $g_n\in\Cal O^*(U)$ are sequences such that $f_n\rightarrow f$
and $g_n\rightarrow g$, then the limit
$$\lim_{n\rightarrow\infty}\{f_n,g_n\}_D,$$
exists, and it is equal to $\{f,g\}_D$.
\enddefinition
\proclaim{Weil's Reciprocity Law 2.4} Let $f$, $g$ be two non-zero
rational functions on $\Bbb P^1$ defined over the field $\Bbb C$. Then
the product of all tame symbols of the pair $(f,g)$ is equal to 1:
$$\prod_{x\in\Bbb P^1}\{f,g\}_x=1.$$
\endproclaim
\definition{Proof} See [12], Proposition 6, pages 44-46. Although it
holds for smooth projective algebraic curves in general, we will
only use this result in the case when the curve is $\Bbb P^1$, when there
is a simple direct proof as follows. The multiplicative group of the
function field of $\Bbb P^1$ is generated by the elements $c\in\Bbb
C^*$, and $z-a$, $a\in\Bbb C$. Since the tame symbols are bilinear and
alternating, we only have to check the identity in the claim for pairs of
these elements. This reduces our proof to three cases: $(z-a,z-b)$, when
$a\neq b$, $(z-a,z-a)$ and $(c,z-a)$. We compute:
$$\split\prod_{x\in\Bbb P^1}\{z-a,z-b\}_x=&\{z-a,z-b\}_a\{z-a,z-b\}_b
\{z-a,z-b\}_{\infty}\\=&(a-b)^{-1}(b-a)(-1)=1,\\
\prod_{x\in\Bbb P^1}\{z-a,z-a\}_x=&\{z-a,z-a\}_a\{z-a,z-a\}_{\infty}
=(-1)(-1)=1,\\
\prod_{x\in\Bbb P^1}\{c,z-a\}_x=&\{c,z-a\}_a\{c,z-a\}_{\infty}
=cc^{-1}=1.\text{\ }\square\endsplit$$
\enddefinition
\definition{Definition 2.5} Let $U$ be a connected rational subdomain of $\Bbb P^1$. For any $D\in\partial U$ and for any two $f$, $g\in\Cal R^*(U)$ we define the rigid analytical regulator $\{f,g\}_D$ by the formula:
$$\{f,g\}_D=\prod_{x\in D}\{f,g\}_x\in\Bbb C^*.$$
Since the product on the right hand side is finite, the regulator is
well-defined. For each $D\in\partial U$ and $f\in\Cal R(U)$ let
$\deg_D(f)$ denote the number of zeros $z$ of $f$ with $z\in D$ counted
with multiplicities minus the number of poles $z$ of $f$ with $z\in D$
counted with multiplicities. For every real number $0<\epsilon<1$ we define the set $\Cal R_{\epsilon}(U)$ as the intersection $\Cal R_{\epsilon}=\Cal O_{\epsilon}(U)\cap\Cal R(U)$. We also define the sets $\Cal O_1(U)=\bigcup_{0<\epsilon<1}\Cal O_{\epsilon}(U)$ and $\Cal R_1(U)=\bigcup_{0<\epsilon<1}\Cal R_{\epsilon}(U)$.
\enddefinition
\proclaim{Lemma 2.6} (i) The regulator is a bilinear map:
$\{\cdot,\cdot\}_D:\Cal R^*(U)\otimes\Cal R^*(U)\rightarrow\Bbb C^*$,

\noindent (ii) the regulator $\{\cdot,\cdot\}_D$ is alternating:
$\{f,g\}_D\cdot\{g,f\}_D=1$,

\noindent (iii) if $f$, $1-f\in\Cal R^*(U)$, then $\{f,1-f\}_D$ is $1$.

\noindent (iv) for any $f\in\Cal R^*(U)$, $c\in\Bbb C^*$ we have
$\{c,f\}_D=c^{\deg_D(f)}$.
\endproclaim
\definition{Proof} The first three claims hold because they hold for the
tame symbol. Claim $(iv)$ is obvious.\ $\square$
\enddefinition
\proclaim{Lemma 2.7} The set $\Cal R_{\epsilon}(U)$ is a subgroup
of $\Cal R^*(U)$.
\endproclaim
\definition{Proof} Note that $\Cal R_{\epsilon}(U)$ is a subset of $\Cal
R^*(U)$, so we have to show that it is a group with respect to  multiplication. If
$f\in\Cal R_{\epsilon}(U)$ then $|f(z)|=1$ for each $z\in U$ by the
ultrametric inequality. Hence
$$\|1-f^{-1}\|=\sup_{z\in U}|f(z)|^{-1}|f(z)-1|\leq
\epsilon.$$
Similarly for any $f$, $g\in\Cal R_{\epsilon}(U)$ we have
$$\split\|1-fg\|\leq&\sup_{z\in U}\max(|1-f(z)|,
|f(z)-f(z)g(z)|)\\\leq&\max(\sup_{z\in U}|1-f(z)|,\sup_{z\in
U}|f(z)|\cdot|1-g(z)|)\leq\epsilon\text{.\ }\square\endsplit$$
\enddefinition
\definition{Definition 2.8} We call a set $D\subset\Bbb P^1$ a closed disk if it is the image of the set $\{z\in\Bbb C||z|\leq1\}$ under a M\"obius transformation. This terminology might be confusing to some as closed disks are open and closed in the natural topology of $\Bbb P^1$, just like open disks. On the other hand closed disks are rational subdomains while open disks are not. For every open disk $D$ let $\overline D$ denote the unique closed disk which contains $D$ and minimal with respect to this property. Moreover let $\partial D$ denote the complement of $D$ in $\overline D$: this notation will not cause confusion because open disks are not rational subdomains. For every $z\in\Bbb C$ and positive number $\rho\in|\Bbb C|$ let $D(z,\rho)$, $\overline D(z,\rho)$, $D(\infty,\rho)$ and $\overline D(\infty,\rho)$ denote the following open disks and closed disks:
$$D(z,\rho)=\{z\in\Bbb C||z|<\rho\}\text{, }\overline D(z,\rho)=\{z\in\Bbb C||z|\leq\rho\}\text{ and }$$
$$D(\infty,\rho)=\{z\in\Bbb C||z|>\rho^{-1}\}\cup\{\infty\}\text{, }\overline D(\infty,\rho)=\{z\in\Bbb C||z|\geq\rho^{-1}\}\cup\{\infty\}.$$
Let $\Bbb C^0$, $\Bbb C^{00}$ and $\bold k$ denote the valuation ring $\{z\in\Bbb C||z|\leq1\}$ of $\Bbb C$, its maximal ideal $\{z\in\Bbb C||z|<1\}$, and its residue field $\Bbb C^0/\Bbb C^{00}$, respectively. Let $\overline z\in\bold k$ denote the reduction of every $z\in\Bbb C^0$ modulo the ideal $\Bbb C^{00}$. For every finite set $S\subset\bold k$ let $U(S)$ denote the set $\{z\in\Bbb C^0|\overline z\notin S\}$. This set is a connected rational subdomain of $\Bbb P^1$, because it is the complement of the pair-wise disjoint disks $D(\infty,1)$ and $D(s,1)$, where $s$ is an element of $\Cal S\subset\Bbb C^0$, a set of representatives of the
residue classes in $S$. Finally let $S^c$ denote the complement of every subset $S$ of $\Bbb P^1$.
\enddefinition
\proclaim{Lemma 2.9} Assume that $0\in\Cal S$. Then every $f\in\Cal
R_{\epsilon}(U(S))$ can be written in the form:
$$f(z)=c\prod_{a\in D(\infty,1)}(1-{z\over a})^{v(a)}\cdot\prod_{s\in\Cal
S}\prod_{a\in D(s,1)}(1-{a-s\over z-s})^{v(a)},$$
where $c\in U_{\epsilon}$ and for each $a\in\Bbb C$ the integer $v(a)$ is
the multiplicity of $a$ in the divisor of $f$.
\endproclaim
\definition{Proof} The rational function:
$$g(z)=\prod_{a\in D(\infty,1)}(1-{z\over a})^{v(a)}\cdot\prod_{s\in\Cal
S}\prod_{a\in D(s,1)}(1-{a-s\over z-s})^{v(a)},$$
is in $\Cal R_1(U(S))$, because it is a product of elements of $\Cal
R_1(U(S))$, which is a group, since it is the union of a chain of groups.
The rational function $f(z)/g(z)$ is also in $\Cal R_1(U(S))$, but it has
no zeros or poles in $\Bbb C-\Cal S$, so it must be equal to the function
$c\prod_{s\in\Cal S}(z-s)^{n(s)}$ for some $c\in\Bbb C^*$ and
$n(s)\in\Bbb Z$. Since $\|f/g\|=1$, we have
$$1=\sup_{z\in U(S)}\left|c\prod_{s\in\Cal S}(z-s)^{n(s)}\right|=|c|.$$
Hence $c\in\Bbb C^0$, so $f/g$ is a rational function with coefficients
on $\Bbb C^0$. Its reduction modulo the ideal $\Bbb C^{00}$
is $r(z)=\overline c\prod_{s\in S}(z-s)^{n(s)}$. Since $f(z)/g(z)$ is in
$\Cal R_1(U(S))$, the rational function $r(z)$ is identically one on
$\bold k-S$. Hence it is constant as $\bold k$ is algebraically closed,
therefore $f(z)=cg(z)$. By the Proposition of I.1.3 in [4], page 7, we
know that every
$h(z)\in\Cal O(U(S))$ has a generalized Laurent expansion:
$$h(z)=\sum_{n=0}^{\infty}a_nz^n+\sum_{s\in\Cal S}\sum_{n=1}^{\infty}
b^s_n(z-a)^{-n},$$
and $\|h(z)\|=\max(\max_{n=0}^{\infty}|a_n|,
\max_{s\in\Cal S}(\max_{n=1}^{\infty}|b^s_n|))$. The constant term $a_0$
in the generalized Laurent expansion of $g(z)$ is 1, so the constant
term  of $f(z)$ is $c$ and we have $|1-c|\leq\|1-f(z)\|\leq
\epsilon$.\ $\square$
\enddefinition
\proclaim{Proposition 2.10} Every $f\in\Cal R_{\epsilon}(U)$ can be
written in the form:
$$f(z)=\prod_{D\in\partial D}f_D(z),$$
where $f_D(z)\in\Cal R_{\epsilon}(D^c)$ for all $D\in\partial U$, and
these functions are uniquely determined up to a constant factor in
$U_{\epsilon}$.
\endproclaim
\definition{Proof} It is clear that the functions $f_D(z)$ are uniquely
determined up to a factor in $U_{\epsilon}$, so we only have to show
that they exist. Assume first that $U$ is of the form $U(S)$ for some
finite set $S\subset\bold k$. The claim holds trivially for the domain
$U(S)$ if the set $S$ is empty. Otherwise we might assume that $0\in\Cal
S$ by a linear change of coordinates. In this case we know that
$$f(z)=f_{\infty}(z)\cdot\prod_{s\in\Cal S}f_s(z),$$
where $f_{\infty}(z)\in\Cal R(\overline D(0,1))\cap\Cal R_1(U)$, and
$f_s(z)\in\Cal R(\overline D(\infty,1))\cap\Cal R_1(U)$ for all $s\in\Cal
S$ by Lemma 2.9. Moreover $f_{\infty}(0)\in U_{\epsilon}$ and
$f_s(\infty)\in U_{\epsilon}$ for all $s\in\Cal S$. However the functions
$f_{\infty}(z)$ and $f_s(z)$ are unique up to a factor in $U_{\epsilon}$,
so $f_{\infty}(a)$ and $f_s(b)\in U_{\epsilon}$ for all $a\in D(0,1)$
and $b\in D(\infty,1)$ for all $s\in\Cal S$, which can be seen by
anticipating the automorphisms $z\mapsto z+a$ and ${1\over z}
\mapsto{1\over z} +{1\over b}$ of $U$, respectively. Let
$1-f_{\infty}(z)=\sum_{n=0}^{\infty}a_n z^n$. There is an
$N(\epsilon)\in\Bbb N$ such that $|a_n|<\epsilon$ if
$n>N(\epsilon)$. Let $B=\{b\in D(0,1)||b|^{n-m}\neq|a_m/a_n|,
\forall n,\forall m\leq N(\epsilon)\}$. If $b\in B$, then
$$\max_{n\leq N(\epsilon)}|a_n|\cdot|b|^n=|\sum_{n\leq N(\epsilon)}a_n
b^n|\leq\max(|1-f_{\infty}(b)|,|\sum_{n>N(\epsilon)}a_nb^n|)
\leq\epsilon.$$
Since $\sup_{b\in B}|b|=1$, the inequality above implies that $\max_{n
\leq N(\epsilon)}|a_n|\leq\epsilon$, so $\|1-f_{\infty}(z)\|
=\max_{n\in\Bbb N}|a_n|\leq\epsilon$. Hence $f_{\infty}\in
\Cal R_{\epsilon}(\overline D(0,1))$, and a similar argument shows that
$f_s\in\Cal R_{\epsilon}(\overline D(\infty,1))$ for all $s\in\Cal S$, so the
claim holds for domains of the form $U(S)$.

In the general case we prove the proposition by induction on the
cardinality of $\partial U$. When $\partial U$ is empty the claim is
obvious. Otherwise let $D$ be a boundary component of $U$. If $\partial
D\nsubseteqq U$ then there is a $D'\in\partial U$ such that $\emptyset
\neq\partial D\cap D'$. We claim that $D'\subset\partial D$ whenever $D$, $D'$ are two disjoint disks such that $\emptyset\neq\partial D\cap D'$, and $D\cup D'\neq\Bbb P^1$. We may assume that $\infty\notin D\cup D'$ by a
linear change of coordinates. Let $a\in\partial D\cap D'$: then $\overline
D=\overline D(a,\rho_1)$ and $D'=D(a,\rho_2)$ for some $\rho_1$,
$\rho_2$, where $\overline D$ denotes the closure of $D$. If
$\rho_2>\rho_1$, then $D\subset\overline D(a,\rho_1)\subset
D(a,\rho_2)$ which is impossible. Therefore $\rho_2\leq\rho_1$, so
$D'=D(a,\rho_2)\subseteq\overline D(a,\rho_1)-D=\partial D$. We
define the relation $D'\leq D$ on the set $\partial U$ by the rule
$D'\subset\partial D$. This is clearly a partial ordering. Let $\Cal D
\subseteq\partial U$ be an equivalence class of minimal elements respect
to this ordering. This means that all $D\in\Cal D$ is minimal with respect to 
this ordering, and any $D'\in\partial U$ is an element of $\Cal D$ if and
only if $D\leq D'$ and $D'\leq D$ for any (and hence all) $D\in\Cal D$.

Choose a boundary component $D\in\Cal D$ and make a linear change of
coordinates such that $D=D(0,1)$. Then all $D'\in\Cal D$ is of the
form $D(s,1)$, where $s\in\Cal S$, and the latter is a finite subset of
$\Bbb C^0$. Since these disks are pair-wise disjoint, $\Cal S$ injects
into $\bold k$ with respect to  the reduction modulo the ideal $\Bbb C^{00}$.
If $S$ denotes the image of $\Cal S$ with respect to  this map, then the set
$E=\bigcap_{D\in\Cal D}\partial D$ is equal to $U(S)$. Therefore $E$ is an
affinoid subdomain of $U$, and we have already proved in the previous
paragraph that the claim holds for this domain. Hence
$$f(z)=g_{\Cal D}(z)\cdot\prod_{D\in\Cal D}f_D(z)$$
on this set, where $g_{\Cal D}(z)\in\Cal R_{\epsilon}(\overline
D(0,1))$ and $f_D(z)\in\Cal R_{\epsilon}(D^c)$ for all $D\in\Cal D$.
Since the set $E$ is infinite, the equation above holds for all $z\in\Bbb
P^1$. Hence
$$g_{\Cal D}(z)={f(z)\over\prod_{D\in\Cal D}f_D(z)}\in\Cal
R_{\epsilon}(U),$$
so $g_D(z)\in\Cal R_{\epsilon}(Y)$, too, where $Y=U\cup\overline
D(0,1)$. Because $|\partial Y|<|\partial U|$, the induction hypothesis
implies that $g_D(z)=\prod_{D'\in\partial Y}f_{D'}(z)$ with
$f_{D'}(z)\in\Cal R_{\epsilon}(D'{}^c)$. Since $\partial U=\partial
Y\cup\Cal D$, the claim follows.\ $\square$
\enddefinition
\proclaim{Proposition 2.11} For every boundary component $D\in\partial U$
and for each $f\in\Cal R_{\epsilon}(U)$ and $g\in\Cal R^*(U)$
we have $\{f,g\}_D\in U_{\epsilon}$.
\endproclaim
\definition{Proof} If the cardinality $|\partial U|=1$ then the
claim is clear, since by Weil's reciprocity $\{f,g\}_D=1$. Hence we may
assume that $\infty\in D^c-U$ by choosing an appropriate coordinate
function $z$. Since $U_{\epsilon}$ is a group, it is sufficient to
prove the claim for $g(z)$ equal to one of the following generators of
of $\Cal R^*(U)$: $c\in\Bbb C^*$, $z-a$, where $a\in D$, and $z-b$, where
$b\in D^c-U$. Then
$$\{f(z),c\}_D=\prod_{D'\in\partial U}\{f_{D'}(z),c\}_D=\{f_D(z),c\}_D=
c^{-\deg_D(f_D)}.$$
But $\deg_D(f_D)$ is zero, since $f_D$ has no zeros or poles in $D^c$,
hence $\{f(z),c\}_D=1$. Also
$$\split\{f(z),z-a\}_D=&\{f_D(z),z-a\}_D\cdot\prod_{D'\in\partial U-\{D\}}
\{f_{D'}(z),z-a\}_D\\=&\prod_{x\in D^c}\{f_D(z),z-a\}^{-1}_x\cdot
\prod_{D'\in\partial U-\{D\}}\{f_{D'}(z),z-a\}_a\\
=&\{f_D(z),z-a\}_{\infty}^{-1}\cdot\!\!\!\!\!\!\prod_{D'\in\partial
U-\{D\}}\!\!\!\!\!\!f_{D'}(a)=f_D(\infty)\cdot\!\!\!\!\!\!\prod_{D'
\in\partial U-\{D\}}\!\!\!\!\!\!f_{D'}(a)\in U_{\epsilon},\endsplit$$
where we used Weil's reciprocity law in the second line. Finally
$$\split\{f(z),z-b\}_D=&\prod_{D'\in\partial
U}\{f_{D'}(z),z-b\}_D=\{f_D(z),z-b\}_D=\prod_{x\in D^c}\{f_D(z),
z-b\}^{-1}_x\\=&\{f_D(z),z-b\}^{-1}_b\cdot\{f_D(z),z-b\}_{\infty}^{-1}=
f_D(b)^{-1}\cdot f_D(\infty)\in U_{\epsilon},\endsplit$$
using again Weil's reciprocity law.\ $\square$
\enddefinition
\definition{Proof of Theorem 2.2} Let $f$ and $g$ be elements of $\Cal O^*(U)$
and let $D\in\partial U$ be a boundary component. By definition (see [4],
page 5) there are sequences $f_n\in\Cal R^*(U)$, $g_n\in\Cal R^*(U)$ such
that $f_n\rightarrow f$ and $g_n\rightarrow g$. We say that the rigid
analytical regulator $\{f,g\}_D\in\Bbb C^*$ exists, if the limit
$$\lim_{n\rightarrow\infty}\{f_n,g_n\}_D,$$
exists, it is an element of $\Bbb C^*$, and it is independent of the
choice of the sequences $f_n$ and $g_n$. In this case we define
$\{f,g\}_D$ to be this limit. We start our proof by showing that the rigid analytical $\{f,g\}_D$ regulator exists. Let $f_n\in\Cal R^*(U)$ and $g_n\in\Cal R^*(U)$ be two sequences such that $f_n\rightarrow f$ and $g_n\rightarrow g$. We will first show that the sequence $\{f_n,g_n\}_D$ converges using Cauchy's criterion. For every $0<\epsilon<1$ there is an
$n(\epsilon)\in\Bbb N$ such that $f_n/f_m\in\Cal R_{\epsilon}(U)$,
$g_n/g_m\in\Cal R_{\epsilon}(U)$ for any $n$, $m\geq n(\epsilon)(U)$. By
Proposition 2.11 for any $n$, $m\geq n(\epsilon)$:
$$\{f_n,g_n\}_D/\{f_m,g_m\}_D=\{f_n/f_m,g_n\}_D\cdot\{f_m,g_n/g_m\}_D\in
U_{\epsilon},$$
so $\{f_n,g_n\}_D$ is a Cauchy sequence. If $f'_n\in\Cal R^*(U)$ and
$g'_n\in\Cal R^*(U)$ are another two sequences such that $f'_n\rightarrow
f$ and $g'_n\rightarrow g$, then the sequences $f''_n$, $g''_n$ also
converge to $f$ and $g$ respectively, where $f''_{2n+1}=f_n$,
$f''_{2n}=f'_n$, and $g''_{2n+1}=g_n$, $g''_{2n}=g'_n$. The limit
$\lim_{n\rightarrow\infty}\{f''_n,g''_n\}_D$ also exists, so it must be
equal to the limit of its subsequences:
$$\lim_{n\rightarrow\infty}\{f''_n,g''_n\}_D=\lim_{n\rightarrow\infty}
\{f_n,g_n\}_D=\lim_{n\rightarrow\infty}\{f'_n,g'_n\}_D,$$
hence the limit is independent of the sequences chosen. On the other hand
$$1=\lim_{n\rightarrow\infty}\{f_n,g_n\}_D\cdot\{g_n,f_n\}_D=
\lim_{n\rightarrow\infty}\{f_n,g_n\}_D\cdot\lim_{n\rightarrow\infty}
\{g_n,f_n\}_D,$$
so this limit is non-zero, so the existence is proved.

Incidentally we also proved that the rigid analytical regulator
satisfies property $(iii)$. Because $\Cal O^*(U)$ is open in $\Cal O(U)$, and
the inverse map is continuous, properties $(ii)$, $(iv)$ and $(v)$ follow
from Lemma 2.6 and Proposition 2.11, by continuity. Property $(i)$ follows
by taking the sequences $f_n=f$ and $g_n=g$. On the other hand
properties $(ii)$ and $(v)$ imply that any map satisfying these
properties must be continuous in both variables, so it is equal to the
rigid analytical regulator if it also satisfies $(i)$.\ $\square$
\enddefinition

\heading 3. The invariance theorem
\endheading

\definition{Definition 3.1} In this chapter we will continue to use the notation of the previous chapter. Let $U$ be a connected rational subdomain
of $\Bbb P^1$, and $f$, $g$ are two meromorphic functions on $U$. Then
 for all $x\in U$ the functions $f$ and $g$ have Laurent series expansion
around $x$, in particular their tame symbol $\{f,g\}_x$ at $x$ is
defined. Let $\Cal M(U)$ denote the field of meromorphic functions of
$U$. For every $x\in U$ the tame symbol at $x$ extends to a homomorphism
$\{\cdot,\cdot\}_x:K_2(\Cal M(U))\rightarrow\Bbb C^*$. We define the
group $K_2(U)$ as the kernel of the direct sum of tame symbols:
$$\bigoplus_{x\in U}\{\cdot,\cdot\}_x:K_2(\Cal M(U))\rightarrow
\bigoplus_{x\in U}\Bbb C^*.$$
Let $k=\sum_if_i\otimes g_i\in K_2(\Cal M(U))$, where $f_i$, $g_i\in\Cal M(U)$, and let $D\in\partial U$. Let moreover $Y$ be a connected rational
subdomain of $U$ such that $f_i$, $g_i\in\Cal O^*(Y)$ for all $i$ and
$\partial U\subseteq\partial Y$. Define the rigid analytical regulator
$\{k\}_D$ by the formula:
$$\{k\}_D=\prod_i\{f_i|_Y,g_i|_Y\}_D.$$
\enddefinition
\proclaim{Theorem 3.2} $(i)$ For each $k\in K_2(\Cal M(U))$ the rigid analytical regulator $\{k\}_D$ is well-defined, and it is a homomorphism
$\{\cdot\}_D:K_2(\Cal M(U))\rightarrow\Bbb C^*$,

\noindent $(ii)$ for any two functions $f$, $g\in\Cal O^*(U)$ we have
$\{f\otimes g\}_D=\{f,g\}_D$,

\noindent $(iii)$ for every $k\in K_2(U)$ the product of all regulators
on the boundary components of $U$ is equal to 1:
$$\prod_{D\in\partial U}\{k\}_D=1,$$

\noindent $(iv)$ for every connected subdomain $Y\subseteq U$, boundary component $D\in\partial Y\cap\partial U$ and $k\in K_2(\Cal M(U))$ we have:
$$\{k|_Y\}_D=\{k\}_D.$$
\endproclaim
In order to prove this theorem, we will need three lemmas.
\proclaim{Lemma 3.3} Let $f$, $g\in\Cal O^*(U)$, where $U$ is a
connected rational subdomain of $\Bbb P^1$. Then the product of all rigid
analytical regulators of the pair $(f,g)$ on the boundary components of
$U$ is equal to 1:
$$\prod_{D\in\partial U}\{f,g\}_D=1.$$
\endproclaim
\definition{Proof} Let $f_n$, $g_n$ be two sequences of rational functions
invertible on  $U$ which converge to $f$ and $g$ on the domain,
respectively. Then
$$\split\prod_{D\in\partial U}\{f,g\}_D=&\prod_{D\in\partial U}
\lim_{n\rightarrow\infty}\{f_n,g_n\}_D=\lim_{n\rightarrow\infty}
\prod_{D\in\partial U}\prod_{x\in D}\{f_n,g_n\}_x\\=&
\lim_{n\rightarrow\infty}\prod_{x\notin U}\{f_n,g_n\}_x=1\endsplit$$
by Weil's reciprocity law.\ $\square$
\enddefinition
\proclaim{Lemma 3.4} Let $Y\subseteq U$ be two connected rational subdomains. Then for every $f$, $g\in\Cal O^*(U)$ and $D\in\partial Y$ we have:
$$\{f|_Y,g|_Y\}_D=\prod\Sb D'\in\partial U\\D'\subseteq D\endSb
\{f,g\}_{D'},$$
where the empty product is understood to be equal to 1.
\endproclaim
\definition{Proof} We may reduce immediately to the case when $f$, $g\in\Cal R^*(U)$ by approximation. Then we have:
$$\split\{f|_Y,g|_Y\}_D=&\prod_{x\in D}\{f|_Y,g|_Y\}_x=
\prod_{x\in D}\{f,g\}_x\\=&\prod\Sb D'\in\partial U\\D'\subseteq D\endSb
\left(\prod_{x\in D'}\{f,g\}_x\right)=\prod\Sb D'\in\partial U\\D'\subseteq D\endSb\{f,g\}_{D'}\text{.\ }\square\endsplit$$
\enddefinition
\proclaim{Lemma 3.5} Let $U$ be a connected rational subdomain and let $D\subset U$ be an open disk. Let $f$, $g$ be two meromorphic functions on $U$ which are holomorphic and invertible on the rational subdomain $Y=U-D$. Then we have:
$$\{f|_Y,g|_Y\}_D=\prod_{x\in D}\{f,g\}_x.$$
\endproclaim
\definition{Proof} The meromorphic functions $f$ and $g$ have only finitely many poles and zeros on the disk $D$. Let $S$ denote the set of these points. For every such a pole or zero $s\in S$ choose a small open disk $D_s\subset D$ containing $s$ such that the pair-wise intersections of these disks are empty. Let $W$ denote the rational subdomain $U-\cup_{s\in S}D_s$. Then
$$\{f|_Y,g|_Y\}_D=\prod_{s\in S}\{f|_W,g|_W\}_{D_s}$$
by Lemma 3.4. Hence we may assume without loosing generality that
$D=D(a,\rho)$ and $a$ is the only pole or zero of $f$ and $g$ on $D$ by letting $W$, $D_s$ play the role of $U$ and $D$, respectively, in the claim above. Write $f(z)=(z-a)^kf_0(z)$, $g(z)=(z-a)^lg_0(z)$, where $f_0(z)$, $g_0(z)$ are elements of $\Cal O^*(U)$. Let $f_n$, $g_n$ be two sequences of rational functions invertible on the rational subdomain $U$ which converge to $f_0$, $g_0$, respectively, with respect to the supremum norm on $U$. Then
$$\split\{f,g\}_D=&\lim_{n\rightarrow\infty}\{(z-a)^kf_n,(z-a)^lg_n\}_D=
\lim_{n\rightarrow\infty}\{(z-a)^kf_n,(z-a)^lg_n\}_a\\
=&\{(z-a)^k,(z-a)^l\}_a\cdot\lim_{n\rightarrow\infty}\{(z-a)^k,g_n\}_a
\cdot\lim_{n\rightarrow\infty}\{f_n,(z-a)^l\}_a\\=&
(-1)^{kl}\cdot\lim_{n\rightarrow\infty}g_n(a)^{-k}\cdot
\lim_{n\rightarrow\infty}f_n(a)^l=(-1)^{kl}g_0(a)^{-k}f_0(a)^l=
\{f,g\}_a.\text{\ }\square\endsplit$$
\enddefinition
\definition{Proof of Theorem 3.2} In order to show that the rigid analytical regulator is well-defined, we have to prove that:
\roster
\item"$(a)$" for each $k\in K_2(U)$ there is an affinoid subdomain
$Y$ with the properties required in Definition 3.1,
\item"$(b)$" the value of the rigid analytical regulator is independent
of the choice of the affinoid subdomain $Y$,
\item"$(c)$" the value of the rigid analytical regulator is independent
of the choice of the presentation $k=\sum_if_i\otimes g_i$.
\endroster
Let $S\subset U$ a finite set such that none of the functions $f_i$ and
$g_i$ has a zero or a pole on $U-S$. For each $s\in S$ take an open disk
$D_s\subset U$ containing $s$ such that these disks are pair-wise
disjoint. The set $Y=U-\bigcup_{s\in S}D_s$ is a connected rational
subdomain which satisfies the properties required in Definition 3.1, so claim $(a)$ holds. Let $Y'$ be another such subdomain. It is clear that the set $Y\cup Y'$ is also a connected rational subdomain which satisfies these properties. Lemma 3.4 applied to the inclusions $Y\hookrightarrow Y\cup Y'$ and $Y'\hookrightarrow Y\cup Y'$ implies that $(b)$ is also true. Let $k=
\sum_if_i'\otimes g'_i$ be another presentation of $k$. By definition
$$\sum_if_i\otimes g_i-\sum_if_i'\otimes g'_i=\sum_jr_j$$
in the free group generated by symbols $f\otimes g$, where $f$, $g\in\Cal
M(U)$, and the elements $r_j$ are defining relations of the group
$K_2(\Cal M(U))$, i.e. they are of the form:
$$fg\otimes h-f\otimes h-g\otimes h\text{,\ }h\otimes fg-h\otimes f
-h\otimes g\text{\ ,\ or\ }f\otimes(1-f)\text{,\ }f\neq1\text{.}$$
It is possible to choose a connected rational subdomain $Y\subset U$ which satisfies the properties required in Definition 3.1 and it does not contain any of the zeros or poles of the functions $f_i$, $g_i$, $f_i'$, $g_i'$, and the functions appearing in the relations $r_j$. By $(ii)$ and $(iv)$ of Theorem 2.2 the regulator $\{\cdot,\cdot\}_D$ evaluated on the relations $r_j$ is equal to 1, so the products:
$$\prod_i\{f_i|_Y,g_i|_Y\}_D=\prod_i\{f'_i|_Y,g'_i|_Y\}_D$$
are equal. Therefore $(c)$ holds, too. The same argument also shows that the
map $\{\cdot\}_D:K_2(\Cal M(U))\rightarrow\Bbb C^*$ is a homomorphism. Claim
$(ii)$ is obvious, if we choose $Y=U$ in the definition. We start the proof of claim $(iii)$ by noting that
$$\prod_{D\in\partial Y}\{k|_Y\}_D=\prod_{D\in\partial U}\{k\}_D\prod_{D\in
\partial Y-\partial U}\{k|_Y\}_D$$
by Lemma 3.4, where $Y$ is a connected rational subdomain which satisfies the properties required in Definition 3.1 with respect to some presentation of $k$. By Lemma 3.5 the factors of the second product on the right hand side are all equal to $1$. On the other hand Lemma 3.3 applied to $Y$ implies that the product on the left hand side is equal to $1$. Hence the first product on the right hand side is equal to $1$, too. We may assume that $k=f\otimes g$ for some $f$, $g\in\Cal M(U)$ by bilinearity while we prove claim $(iv)$. There is a connected subdomain $Z\subseteq U$ such that $\partial U\subseteq\partial Z$, the intersection $Y\cap Z$ is nonempty and $f$, $g\in\Cal O^*(Z)$. Clearly $D\in\partial(Y\cap Z)$, so it will be sufficient to prove that
$$\{f|_{Y\cap Z},g|_{Y\cap Z}\}_D=\{f|_Z,g|_Z\}_D.$$
We may even assume that $f$ and $g$ are in $\Cal R^*(Y)$ by approximation. But in this case the claim is obviously true.\ $\square$
\enddefinition
\proclaim{Theorem 3.6} There is a unique set of homomorphisms $\deg_D:\Cal M^*(U)\rightarrow\Bbb Z$ where $U$ is any connected rational subdomain and $D\in\partial U$ is a boundary component with the following properties:
\roster
\item"$(i)$" the homomorphism $\deg_D$ is zero on $\Cal O_1(U)$,
\item"$(ii)$" for every $f\in\Cal R^*(U)$ the integer $\deg_D(f)$ is the quantity defined in 2.5,
\item"$(iii)$" for every $f\in\Cal M^*(U)$ we have $\deg_D(f|_Y)=\deg_D(f)$ where $Y\subseteq U$ is any connected rational subdomain satisfying the property $D\in\partial Y$.
\endroster
\endproclaim
This map is equal to the integer-valued function introduced in [4], Proposition I.3.1, page 19, in the limited context when the latter is defined.
\definition{Proof} By condition $(i)$ it is clear that the homomorphism $\deg_D$ restricted to $\Cal O^*(U)$ must be continuous with respect to the supremum topology on $\Cal O^*(U)$ and the discrete topology on $\Bbb Z$, if it exists. For every $f\in\Cal M^*(U)$ there is a $Y\subseteq U$ connected rational subdomain such that $D\in\partial Y$ and $f|_Y\in\Cal O^*(Y)$, so the set of homomorphisms $\deg_D:\Cal M^*(U)\rightarrow\Bbb Z$ should be unique. Pick an element $c\in\Bbb C$ with $|c|>1$ and define $\deg_D(f)$ as
$$\deg_D(f)=\log_{|c|}(|\{c,f\}_D|)$$
where $\log_a$ is the logarithm with base $a$ for any positive real number $a$. This homomorphism is well-defined by Theorem 3.2. This function also satisfies condition $(ii)$ of the proclaim above by claim $(iv)$ of Lemma 2.3. Hence the image of $\Cal R^*(Y)$ with respect to this map lies in $\Bbb Z$. But it is also dense in the image of $\Cal O^*(Y)$ via this map, so the latter lies in $\Bbb Z$, too. We may conclude that the homomorphisms $\deg_D$ take integral values. By $(v)$ of Theorem 2.2 the element $\{c,f\}_D$ is a unit in $\Bbb C^0$ when $f\in\Cal O_1(U)$, so condition $(i)$ is also satisfied. Condition $(iii)$ is the consequence of claim $(iv)$ of Theorem 3.2.\ $\square$
\enddefinition
\definition{Definition 3.7} For every $U\subset\Bbb P^1$ connected rational subdomain let $\Bbb Z\partial U$ denote the free abelian group with the elements of $\partial U$ as free generators. Let $H_1(U)$ denote the quotient of $\Bbb Z\partial U$ by the $\Bbb Z$-module generated by $\sum_{D\in\partial U}D$. For every $D\in\partial U$ we let $D$ denote the class of $D$ in $H_1(U)$ as well. Let $\Cal Ab$ denote the category of abelian groups. Let $\Cal Crs$ denote the category whose objects are connected rational subdomains of $\Bbb P^1$ and whose morphisms are holomorphic maps between them. Finally for every pair $a\leq b$ of numbers in $|\Bbb C|$ let $A(a,b)$ denote the closed annulus $\Bbb P^1-D(0,a)-D(\infty,1/b)$. Of course it is a connected rational subdomain.
\enddefinition
\proclaim{Theorem 3.8} There is a unique functor $H_1:\Cal Crs\rightarrow\Cal Ab$ with the following properties:
\roster
\item"$(i)$" for every $U\subset\Bbb P^1$ connected rational subdomain $H_1(U)$ is the group defined in 3.7,
\item"$(ii)$" for every map $U\rightarrow Y$ which is the restriction of a projective linear transformation $f$ and $D\in\partial U$ boundary component we have:
$$H_1(f)(D)=f(D)\in H_1(Y),$$
\item"$(iii)$" for every $f:U\rightarrow D(a,b)$ holomorphic map and $D\in\partial U$ boundary component we have:
$$H_1(f)(D)=\deg_D(f)D(0,a)\in H_1(A(a,b)).$$
\endroster
\endproclaim
\definition{Proof} We are going to prove first that this functor is unique. Let $h:U\rightarrow Y$ be a holomorphic map between two connected rational subdomains. We need to show that $H_1(h)$ is uniquely determined by the conditions above. We may assume that $Y$ has at least two boundary components. Fix a boundary component $F\in\partial Y$. Then for every other boundary component $F\neq E\in\partial Y$ there is a projective linear transformation $j_E$ of $\Bbb P^1$ such that $j_E(E)=D(0,1)$ and $\infty\in j_E(F)$. Then $j_E\circ H$ maps into $A(a,b)$ for some $a$ and $b$ for every $F\neq E\in\partial Y$. By property $(ii)$ and $(iii)$ we have:
$$H_1(h)(D)=\sum_{F\neq E\in\partial Y}\deg_D(j_E\circ h)E\in H_1(Y)$$
for every boundary component $D\in\partial U$. In particular this class is uniquely determined. Let $H^1(U)$ denote the quotient
$$H^1(U)=\Cal O^*(U)/(\Bbb C^*\Cal O_1(U)),$$
where $\Bbb C^*\subset\Cal O^*(U)$ is the subgroup of constant functions. For every $f\in\Cal O^*(U)$ let the same letter denote its class in $H^1(U)$ as well. The degree map of Theorem 3.6 induces a bilinear pairing:
$$\deg:H^1(U)\times H_1(U)\rightarrow\Bbb Z$$
characterized by the property:
$$\deg(f,D)=\deg_D(f)$$
for every $f\in\Cal O^*(U)$ and $D\in\partial U$.
\enddefinition
\proclaim{Proposition 3.9} The pairing $\deg$ is perfect.
\endproclaim
\definition{Proof} We need to show the following two claims in order to prove the proposition:
\roster
\item"(i)" for every $f\in\Cal O^*(U)$ if $\deg_D(f)=0$ for every $D\in\partial U$ then $f\in\Bbb C^*\Cal O_1(U)$,
\item"(ii)" for every homomorphism $l:H_1(U)\rightarrow\Bbb Z$ there is a function $f\in\Cal O^*(U)$ such that $l(h)=\deg(f,h)$ for every $h\in H_1(U)$.
\endroster
We may assume that $f\in\Cal R^*(U)$ by approximation while we show the first claim. Let $\sum_{x\in\Bbb P^1}n(x)x$ be the divisor of $f$. Let $g(z)$ be the product $\prod_{D\in\partial D}g_D(z)$, where
$$g_D(z)=\prod_{x\in D}(z-x)^{n(x)}$$
if $\infty\notin D$, and
$$g_D(z)=\prod_{x\in D}({1\over z}-{1\over x})^{n(x)},$$
otherwise. The rational functions $f$ and $g$ have the same divisor, so their quotient is constant. Therefore it will be enough to show that $g\in\Cal O_1(U)$. We will prove that the functions $g_D\in\Cal R_1(U)$ which is sufficient by Lemma 2.7. First consider the case when $\infty\notin D$. In this case $D=D(c,d)$ for some $c\in\Bbb C$ and $d\in|\Bbb C|$. By assumption $\sum_{x\in D}n(x)=\deg_D(f)=0$, so the function $g_D(z)$ is the product of factors of the form $(z-a)/(z-b)$, where $a$, $b\in D$. For any $z\notin D$ we have
$$\left|{z-a\over z-b}\right|={|z-c+c-a|\over|z-c+c-b|}={|z-c|\over|z-c|}=1,$$
so $g_D\in\Cal R_1(U)$ as we claimed. In the other case the argument is similar. We may assume that $\partial U$ has at least two elements while we show the second claim. Fix a boundary component $F\in\partial U$. As $l(F)=-\sum_{F\neq D\in\partial U}l(D)$, there is a rational function $f\in\Cal R^*(U)$ such that $\deg_D(f)=l(D)$ for every other boundary component $F\neq D\in\partial U$. Clearly the function $f$ satisfies the property in the second claim.\ $\square$
\enddefinition
The proof of the existence of the functor $H_1$ is now easy: we define it as the $\Bbb Z$-dual of the contravariant functor $H^1$. Condition $(i)$ is automatic via the identification between $H_1(U)$ and Hom$(H^1(U),\Bbb Z)$ furnished by Proposition 3.9. Properties $(ii)$ and $(iii)$ can be verified by looking at appropriate test functions to compute the effect of $H_1(f)$. In the first case one considers rational functions, in the second case the identity function.\ $\square$
\definition{Definition 3.10} Let $U\subset\Bbb P^1$ be a connected rational subdomain. For every class $c\in H_1(U)$ and element $k\in K_2(U)$ we define the regulator $\{k\}_c$ as
$$\{k\}_c=\prod_{D\in\partial U}\{k\}_D^{c(D)},$$
where $\sum_{D\in\partial U}c(D)D$ is a lift of $c$ in $\Bbb Z\partial U$. By claim $(iii)$ of
Theorem 3.2 this regulator is well-defined. For every holomorphic map $h:U\rightarrow Y$
between two connected rational subdomains let $h^*:K_2(\Cal M(Y))\rightarrow K_2(\Cal M(U))$
be the pull-back homomorphism induced by $h$. By restriction it induces a homomorphism
$K_2(Y)\rightarrow K_2(U)$.
\enddefinition
\proclaim{Theorem 3.11} For any $k\in K_2(Y)$ and $c\in H_1(U)$ we
have:
$$\{h^*(k)\}_c=\{k\}_{H_1(h)(c)}.$$
\endproclaim
\definition{Proof} Let $k=\sum_if_i\otimes g_i\in K_2(\Cal M(Y))$, where $f_i$, $g_i\in\Cal M(U)$. Let moreover $Y'$ be a connected rational
subdomain of $Y$ such that $f_i$, $g_i\in\Cal O^*(Y)$ for all $i$ and
$\partial Y\subseteq\partial Y'$. There is a connected rational
subdomain $U'$ of $U$ such that $h(U')\subseteq Y'$ and
$\partial U\subseteq\partial U'$. The map $H_1(U')\rightarrow H_1(U)$ induced by the inclusion is surjective, so there is a $c'\in H_1(Y')$ whose image is $c$. We claim that $\{h^*(k)|_{U'}\}_{c'}=\{h^*(k)\}_c$. We may write $c'$ as a sum $c'=c_1+c_2$ where $c_1$, $c_2$ can be represented as the linear combination of boundary components lying in $\partial U$ and $\partial U'-\partial U$, respectively. We have $\{h^*(k)|_{U'}\}_{c_1}=\{h^*(k)\}_c$ by definition, on the other hand $\{h^*(k)|_{U'}\}_{c_2}=1$ as $h^*(k)$ is an element of $K_2(U)$. The same argument shows that $\{k|_{Y'}\}_{H^1(h|_{U'})(c')}=\{k\}_{H^1(h)(c)}$, so it will be sufficient to prove the claim for $U'$, $Y'$, $h|_{U'}$, $c'$ and $k|_{Y'}$ instead of $U$, $Y$, $h$, $c$ and $k$, respectively. In other words we may assume that $k=f\otimes g$ for some $f$, $g\in\Cal O^*(Y)$. Let $f_n\in\Cal R^*(Y)$ and $g_n\in\Cal R^*(Y)$ be two sequences such that $f_n\rightarrow f$ and
$g_n\rightarrow g$. Obviously $\|f\circ h-f_n\circ h\|\leq\|f-f_n\|$ for all $n\in\Bbb N$, so $f_n\circ h\rightarrow f\circ h$. The same holds for $g$, so $$\{h^*(f\otimes g)\}_c=\lim_{n\rightarrow\infty}\{h^*(f_n\otimes g_n)\}_c
\text{ and }\{f\otimes g\}_{H_1(h)(c)}=\lim_{n\rightarrow\infty}\{f_n\otimes g_n\}_{H_1(h)(c)}$$
by Remark 2.3. Therefore it is sufficient to show the claim when $f$, $g\in\Cal R^*(Y)$. We may also assume that $\infty\notin U$ and $\infty\notin Y$ by shrinking $Y$ and $U$ the same way as above if necessary.
\enddefinition
\proclaim{Lemma 3.12} For every $f\in\Cal O(U)$ the following holds: for every $\epsilon>0$ there is a $\delta>0$ such that $|f(x)-f(y)|<\epsilon$ for every $x$, $y\in U$ with $|x-y|<\delta$.
\endproclaim
\definition{Proof} Of course the claim above is just the analogue of the usual uniform continuity property. The reason that it is not completely obvious in this case is that $\Bbb C$ is not locally compact. Let $\Cal U(U)$ denote the set of all functions $f\in\Cal O(U)$ which satisfy the property in the claim above. It is clear that $\Cal U(U)$ is a $\Bbb C$-subalgebra of $\Cal O(U)$. Moreover for every $f\in\Cal O^*(U)$ we have the estimate:
$$\left|{1\over f(x)}-{1\over f(y)}\right|=\left|{f(y)-f(x)\over f(x)f(y)}\right|\leq\|f^{-1}\|^2|f(x)-f(y)|,$$
so for every $f\in\Cal O^*(U)\cap\Cal U(U)$ we have $f^{-1}\in\Cal U(U)$, too. Obviously $z-c\in\Cal U(U)$ for every $c\in\Bbb C$, so $\Cal R(U)\subseteq\Cal U(U)$ by the above. On the other hand $\Cal U(U)$ is closed with respect to the supremum topology, so it must be equal to the whole algebra $\Cal O(U)$.\ $\square$
\enddefinition
Let us return to the proof of Theorem 3.11. Since $\infty\notin U$ there is a rational $\epsilon>0$ such that for every $x\in U$ the disk $D(x,\epsilon)\subset U$. Hence we may choose an infinite sequence $h_n\in \Cal R(U)$ converging to $h$ in the supremum topology such that $h_n(U)\subseteq Y$ for all $n$. By the lemma above $f\circ h_n\rightarrow f\circ h$ and $g\circ h_n\rightarrow g\circ h$ in the supremum topology. Therefore it is sufficient to prove theorem when $h\in\Cal R(U)$, too. We may also assume that $c=C$ for some boundary component $C\in\partial U$ by linearity. Let $F\in\partial Y$ be the unique boundary component which contains $\infty$. We may write $H^1(h)(C)$ uniquely in the form:
$$H^1(h)(C)=\sum_{F\neq D\in\partial Y}n(D)D$$
for some $n(D)\in\Bbb Z$. There is a closed annulus $A(a,b)$ such that $h(z)-c$ maps $U$ into $A(a,b)$ for every $F\neq D\in\partial Y$ and $c\in D$.
Fix a boundary component $F\neq D\in\partial Y$ and for every $c\in D$ let $z(c)$ denote the number of zeros of the rational function $h(z)-c$ lying in the open disk $C$ counted with multiplicities. We claim that $z(c)$ is independent of the choice of $c$. First note that the number of poles of the rational function $h(z)-c$ lying in the open disk $C$ counted with multiplicities is independent of the choice of $c$. This number does not even depend on $D$, and it will be denoted by $p(C)$. By claim $(iii)$ of Theorem 3.8 we have $z(c)=\deg_C(h(z)-c)+p(D)=n(D)+p(C)$ which is clearly independent of the choice of $c$. Let $z(D)$ denote this number. We also claim that for every $\infty\neq c\in F$ the number $z(c)$ of zeros of the rational function $h(z)-c$ lying in the open disk $C$ counted with multiplicities is equal to $p(C)$. We have $F=D(\infty,d)$ for some rational number $d>0$. Hence $\|h(z)\|\leq1/d$, so we have
$$\left\|1-{h(z)-c\over c}\right\|=\left\|{h(z)\over c}\right\|<1,$$
so $\deg_C(h(z)-c)=0$. But $z(c)=\deg_D(h(z)-C)+p(C)=p(C)$ by definition.
For all $x\in\Bbb P^1$ let $v(x)\in\Bbb N$ denote the degree of ramification of the map $h:\Bbb P^1\rightarrow\Bbb P^1$ at $x$. Then
$$\{f\circ h,g\circ h\}_x=\{f,g\}_{h(x)}^{v(x)}$$
for all $x\in\Bbb P^1$. Therefore
$$\split\{h^*(f\otimes g)\}_C=&\prod_{x\in C}\{f\circ h,g\circ h\}_x\\
=&\prod_{x\in C}\{f,g\}^{v(x)}_{h(x)}\\=&
\prod_{D\in\partial Y}\prod_{y\in D}\prod\Sb x\in C\\ h(x)=y\endSb
\{f,g\}^{v(x)}_y\\
=&\prod_{F\neq D\in\partial Y}\prod_{y\in D}\{f,g\}_y^{z(D)}\cdot
\prod_{y\in F}\{f,g\}_y^{p(C)}\\
=&\prod_{F\neq D\in\partial Y}\prod_{y\in D}\{f,g\}_y^{z(D)-p(C)}=
\{f\otimes g\}_{H^1(h)(C)},\endsplit$$
where we used Weil's reciprocity law in the fifth equation.\ $\square$
\definition{Remarks 3.13} The theorem above incorporates two fundamental properties of the
rigid analytic regulator which might be called as biholomorphic invariance and homotopy
invariance based on the analogy explained in 1.1. The fist claims that every biholomorphic map
$h:U\rightarrow U$, where $U$ is any connected rational subdomain of $\Bbb P^1$, which induces
the identity map on $H^1(U)$ leaves the rigid-analytic regulator invariant. The second claims
that whenever we have two connected rational subdomains $Y\subseteq U$ and a boundary
component $D\in\partial U$ such that there is a unique boundary component $D'\in\partial Y$
containing $D$ then the rigid analytic regulators taken at $D$ and $D'$ do not differ for
elements of $K_2(U)$. Of course the best way to formulate these properties is the way we did,
expressing them as a functoriality property via the homology group $H_1$. The latter has a more
high-brow definition using the \'etale cohomology of rigid analytic spaces, but for our
purposes our elementary definition was more suitable.
\enddefinition

\heading 4. Relation to the generalized tame symbol
\endheading

\definition{Notation 4.1} Let $\Bbb R$ be a closed subfield of $\Bbb C$: it is automatically complete with respect to $|\cdot|$. Let $U$ be a connected rational subdomain of $\Bbb P^1$ defined over $\Bbb R$. This means that
$$U=\{z\in\Bbb P^1||f_i(z)|\leq1\ (\forall i=1,\ldots,n)\}$$
for some natural number $n$ and rational functions $f_1,\ldots,f_n\in\Bbb R((t))$. Let $\Cal O_{\Bbb R}(U)$, $\Cal R_{\Bbb R}(U)$, $\Cal O_{\Bbb R}^*(U)$, $\Cal R_{\Bbb R}^*(U)$ and $\Cal M_{\Bbb R}(U)$ denote the algebra of holomorphic functions, the subalgebra of restrictions of $\Bbb R$-rational functions, the groups of invertible elements of these algebras and the field of meromorphic functions on the rigid analytic space $U$, respectively. Let $U$ denote the underlying rational subdomain over $\Bbb C$ by slight abuse of notation. An $\Bbb R$-rational boundary component of $U$ is a set $D\in\partial U$ such that $D$ is of the form $D(a,|\rho|)$ or $D(\infty,|\rho|)$ for some
$a$, $\rho\in\Bbb R$. Let $K_2(U)_{\Bbb R}$ denote the largest subgroup of $K_2(\Cal M_{\Bbb R}(U))$ which maps into $K_2(U)$ under the natural homomorphism $K_2(\Cal M_{\Bbb R}(U))\rightarrow K_2(\Cal M_{\Bbb C}(U))$.
\enddefinition
\proclaim{Proposition 4.2} Let $D$ be an $\Bbb R$-rational boundary
component of $U$, and let $k\in K_2(\Cal M_{\Bbb R}(U))$. Then $\{k\}_D\in
\Bbb R^*$.
\endproclaim
\definition{Proof} We may assume that $k=f\otimes g$ for some $f$,
$g\in\Cal O_{\Bbb R}^*(U)$ by linearity. Since $\Cal R_{\Bbb R}^*(U)$ is dense in $\Cal O_{\Bbb R}^*(U)$ and $\Bbb R$ is complete, we may assume that $f$ and $g$ are actually in $\Cal R^*(U)$ by approximation. We may also
assume that $\infty\in D$ after an $\Bbb R$-linear change of
coordinates. Then by bilinearity we may assume that $f$ and $g$ are
irreducible polynomials in $\Bbb R[t]$. Assume first that $f$ and $g$
are separable, too. Clearly $\{f,g\}_x$ can be different from $1$ only if $x=\infty$ or $x$ is a zero of $f$ or $g$. In the latter case $x$ is an element of the separable closure $\overline{\Bbb R}$ hence $\{f,g\}_x\in
\overline{\Bbb R}$, too. Moreover if $x\in D$ then $x^{\gamma}\in D$, too, where $\gamma$ is any element of Gal$(\overline{
\Bbb R}|\Bbb R)$. Also $\{f,g\}_x^{\gamma}=\{f,g\}_{x^{\gamma}}$
for any $x\in\overline{\Bbb R}$ and $\gamma\in\text{Gal}(\overline{
\Bbb R}|\Bbb R)$. Therefore the product $\prod_{x\in D}\{f,g\}_x$
is an element of $\overline{\Bbb R}$ invariant under the action of
$\text{Gal}(\overline{\Bbb R}|\Bbb R)$, so it is in $\Bbb R$. If $f$
is not separable, then $f=(f')^{p^n}$ for some separable polynomial $f'$
whose coefficients are in a purely inseparable extension $\Bbb L$ of
$\Bbb R$ such that $\Bbb L^{p^n}\subseteq\Bbb R$ where $p$ is the characteristic of $\Bbb R$. It is enough to show that $\{f',g\}_D\in\Bbb L^*$, since $\{f,g\}_D=\{f',g\}_D^{p^n}\in\Bbb R^*$ in this case. The latter follows from applying the same argument to $g$ over $\Bbb L$ and what we have just proved above.\ $\square$
\enddefinition
\definition{Definition 4.3} Let $\bold f$, $\Cal O$ and $\goth m$ denote the residue field of $\Bbb R$, the valuation ring of $\Bbb R$ and the maximal proper ideal of $\Cal O$,  respectively. A finite subset $S\subset\bold k$ is called $\bold f$-rational if it is the zero set of a polynomial with coefficients in $\bold f$. If $S\subset\bold k$ is $\bold f$-rational then the set $U(S)$ introduced in Definition 2.8 is a connected rational subdomain defined over $\Bbb R$.
\enddefinition
\proclaim{Lemma 4.4} Let $S\subset\bold k$ be an $\bold f$-rational subset. Then
\roster
\item"(i)" for every $f\in\Cal M_{\Bbb R}(U(S))^*$ there is an $\bold f$-rational subset $S\subseteq S'\subset\bold k$ such that $f|_{U(S')}$ can be written
in the form:
$$f|_{U(S')}=c(f)f_0,$$
where $c(f)\in\Bbb R$ and $f_0\in\Cal O^*(U(S'))$ with $|f_0(z)|=1$ for all $z\in U(S')$,
\item"$(ii)$" the positive number $|c(f)|$ does not depend on the choice of $S'$ or $f_0$,
\item"$(iii)$" the map $|\cdot|:\Cal M_{\Bbb R}(U(S))^*\rightarrow|\Bbb C|$ given by the rule $f\mapsto|c(f)|$ is a non-archimedean absolute value on the field $\Cal M_{\Bbb R}(U(S))$.
\endroster
\endproclaim
\definition{Proof} The set $Z\subset U(S)$ of zeros and poles of $f$ is
finite. Hence the reduction of the elements of $Z$ with respect to $\Bbb C^{00}$ is a finite set, too. Since every finite subset of $\bold k$ is contained in a finite $\bold f$-rational subset, we may assume that $f\in\Cal O^*(U(S))$ by enlarging $S$ if necessary. Let $g\in\Cal O^*(U(S))$ be another function such that $\|1-f/g\|<1$. We claim that it is sufficient to prove the claim $(i)$ for $g$ in order to prove it for $f$. We may assume that $g=c(g)g_0$ where $c(g)\in\Bbb R$ and $g_0\in\Cal O^*(U(S))$ with $|g_0(z)|=1$ for all $z\in U(S)$ by enlarging $S$ if necessary. Write $f$ as $f=c(g)f_0$. Then $\|1-f_0/g_0\|=\|1-f/g\|<1$ so $|f_0(z)|=1$ for all $z\in U(S)$. Hence we may assume that $f\in\Cal R^*(U(S))$ by approximation. Also note that the elements of $\Cal M_{\Bbb R}(U(S))^*$ satisfying claim $(i)$ form a subgroup. Therefore we may assume that $f$ is in fact a polynomial. Then we may write $f(z)$ as:
$$f(z)=c(f)\sum_{n=0}^{N}a_nz^n$$
with $a_n\in\Cal O$ and $\max_{n=0}^{N}|a_n|=1$. There is a finite $\bold f$-rational subset $S'$ of $\bold k$ such that the reduction of the polynomial $f_0=\sum_{n=0}^{N}a_nz^n$ is nowhere zero on the complement of $S'$. Clearly $f_0\in\Cal O^*(U(S'))$ with $|f_0(z)|=1$ for all $z\in U(S')$.

This proves claim $(i)$. Assume that $S\subseteq S''\subset\bold k$ is another finite $\bold f$-rational subset such that  $f|_{U(S'')}=c(f)'f_0'$ where $c(f)'\in\Bbb R$ and $f_0'\in\Cal O^*(U(S'))$ with $|f_0'(z)|=1$ for all $z\in U(S'')$. On the set $U(S'\cup S'')=U(S')\cap U(S'')$ we have $c(f)f_0=c(f)'f_0'$. As $\|f_0\|=\|f_0'\|=1$ on this set we must have $|c(f)|=|c(f)'|$ as claim $(ii)$ says. Clearly the map $v$ is a homomorphism, so we only have to show that $|c(f+g)|\leq\max(|c(f)|,|c(g)|)$ for any $f$, $g\in\Cal O^*(U(S))$ with $f+g\neq0$ in order to prove claim $(iii)$. There is an $\bold f$-rational subset $S\subseteq S'\subset\bold k$ such that $|c(f)|=\|f|_{U(S')}\|$, $|c(g)|=\|g|_{U(S')}\|$ and $|c(f+g)|=\|(f+g)|_{U(S')}\|$. The last claim now follows from the strong triangle inequality for the spectral norm $\|\cdot\|$.\ $\square$
\enddefinition
\definition{Definition 4.5} For any pair $S\subseteq S'\subset\bold k$ of finite
$\bold f$-rational subsets the inclusion $U(S)\rightarrow U(S')$ induces an imbedding
$\Cal M_{\Bbb R}(U(S))\rightarrow\Cal M_{\Bbb R}(U(S'))$. Under these inclusions the fields
$\Cal M_{\Bbb R}(U(S))$ form a injective system: let $\Cal M$ denote the inductive limit of
this system. By part $(ii)$ of the lemma above the absolute value $|\cdot|$ of part $(iii)$ is
well-defined on $\Cal M$ and makes the latter a valued field. The residue field of $\Cal M$ is
equal to the rational function field $\bold f(t)$ where $t$ is the reduction of the identity
map $z$ with respect to the maximal ideal in the valuation ring of $\Cal M$. We are going to
need a mild extension of the tame symbol. Let $F$ be a field equipped with a valuation
$\nu:F^*\rightarrow\Bbb Q$ and let $O$ and $\bold r$ denote its valuation ring and its residue
field, respectively. Let moreover $\overline{\cdot}:O\rightarrow\bold r$ denote the reduction
modulo the maximal ideal of $O$. For any pair of elements $f$, $g\in F^*$ we are going to
define their generalized tame symbol $T(f\otimes g)\in\bold r^*\otimes\Bbb Q$ as follows.
There is an element $\pi\in F^*$ such that $f=f_0\pi^{n(f)}$ and $g=g_0\pi^{n(g)}$ for some
integers $n(f)$, $n(g)\in\Bbb Z$ and elements $f_0$, $g_0\in\Cal O^*$.  We let
$$T(f\otimes g)=(-1)\otimes\nu(f)\nu(g)\cdot\overline{f_0}\otimes\nu(g)
\cdot\overline{g_0}\otimes(-\nu(f))\in\bold r^*\otimes\Bbb Q.$$
(We only included the first factor in the product above in order to resemble the usual
formula, but it is always equal to $1$). One may prove the usual way that this symbol is
well-defined and satisfies the Steinberg relation. In particular it induces a homomorphism
$T:K_2(F)\rightarrow\bold r^*\otimes\Bbb Q$ which depends on the choice of normalization of
the valuation $v$ linearly. Let $\nu:\Cal M^*\rightarrow\Bbb Q$ be a valuation corresponding the
absolute value $|\cdot|$ and let $T:K_2(\Cal M)\rightarrow\bold f(t)^*\otimes\Bbb Q$ denote
the corresponding generalized tame symbol. For any $k\in K_2(\Cal M(U))$ and
$s\in S\cap\bold f$ let $T_s(k)\in\Bbb Q$ denote the value of $T(k)$ with respect to the
unique $\Bbb Q$-linear extension of the normalized valuation at the closed point
$s\in\Bbb P^1_{\bold f}$.
\enddefinition
\proclaim{Proposition 4.6} We have
$$\nu(\{k\}_{D(s,1)})=T_s(k)$$
for every $k\in K_2(\Cal M_{\Bbb R}(U(S)))$.
\endproclaim
\definition{Proof} The linear transformation $z\mapsto z-s$ maps $U(S)$ biholomorphically onto
$U(\cup_{x\in S}x-s)$ and interchanges $\{\cdot\}_{D(s,1)}$ and $T_s$ with $\{\cdot\}_{D(0,1)}$
and $T_0$, respectively. Hence we may assume that $s=0$. By linearity we may assume that
$k=f\otimes g$ for some $f$, $g\in \Cal M_{\Bbb R}(U(S))$. Since for every $\bold f$-rational
finite set $S\subseteq S'\subset\bold k$ we have $\{k\}_{D(0,1)}=\{k|_{U(S')}\}_{D(0,1)}$ and
$T_0(k)=T_0(k|_{U(S')})$ we may assume that $f=c(f)f_0$ and $g=c(g)g_0$ where $c(f)$, $c(g)\in
\Bbb R^*$ and $f_0$, $g_0\in\Cal O^*(U(S))$ with $|f_0(z)|=|g_0(z)|=1$ for every $z\in U(S)$ by
enlarging $S$ if necessary. Assume that $\|1-f\|<1$. In this case $T(f\otimes g)=1$ by definition
and $|\{f\otimes g\}_{D(0,1)}|=1$ by $(v)$ of Theorem 2.2. Hence we may assume that $f$ and $g$
are rational functions by approximation. We may even assume that $f$ and $g$ are polynomials
using the bilinearity of both sides of the equation we want to prove. In this case we may assume
by multiplying $c(f)$ and $c(g)$ by an element of $\Cal O^*$, if necessary, such that
$$f_0(z)=z^{n(f)}(1+f_{00}(z))\text{ and }g_0(z)=z^{n(g)}(1+g_{00}(z)),$$
where $n(f)$, $n(g)\in\Bbb Z$ and $f_{00}(z)$, $g_{00}(z)\in\Cal O[z]$ with $\|f_{00}\|<1$ and
$\|g_{00}\|<1$. Therefore we may assume that $f_0$ and $g_0$ are powers of $z$ by applying the
argument we used above. In this case the claim is obvious.\ $\square$
\enddefinition

\heading 5. The rigid analytical regulator for Mumford curves
\endheading

\definition{Definition 5.1} Let $\goth X$ be a Hausdorff topological
space. For any $R$ commutative group let $\Cal M(\goth X,R)$ denote the set
of $R$-valued finitely additive measures on the open and compact
subsets of $\goth X$. When $\goth X$ is compact let $\Cal M_0(\goth X,R)$ denote the set of
measures of total volume zero, that is the subset of those $\mu\in\Cal M(\goth X,R)$ such
that $\mu(\goth X)=0$. For every abelian topological group $M$ let $\Cal C_0(\goth X,M)$
denote the group of compactly supported continuous functions $f:\goth X\rightarrow M$. If $M$
is discrete then every element of $\Cal C_0(\goth X,M)$ is locally constant. In this case for
every $f\in\Cal C_0(\goth X,M)$ and $\mu\in\Cal M(\goth X,R)$ we define the
modulus $\mu(f)$ of $f$ with respect to $\mu$ as the $\Bbb Z$-submodule of
$R$ generated by the elements $\mu(f^{-1}(g))$, where $0\neq g\in M$. We also
define the integral of $f$ on $\goth X$ with respect to $\mu$ as the sum:
$$\int_{\goth X}f(x)\text{d}\mu(x)=\sum_{g\in M}g\otimes\mu(f^{-1}(g))\in
M\otimes\mu(f).$$
\enddefinition
\proclaim{Lemma 5.2} $(a)$ If $f$, $g\in\Cal C_0(\goth X,M)$, then $f
\times g\in \Cal C_0(\goth X,M\times M)$, $\mu(f)$, $\mu(g)$ and
$\mu(f+g)$ are contained in $\mu(f\times g)$ and
$$\int_{\goth X}f(x)+g(x)\text{\rm d}\mu(x)=\int_{\goth X}f(x)
\text{\rm d}\mu(x) +\int_{\goth X}g(x)\text{\rm d}\mu(x)\text{\ in\
}M\otimes\mu(f\times g).$$
\noindent $(b)$ Let $f\in\Cal C_0(\goth X,M)$ and $m\in N$, where $N$ is
also a $\Bbb Z$-module. Then $m\otimes f$ is in $\Cal C_0(\goth X,N\otimes M)$,
$\mu(m\otimes f)\subseteq\mu(f)$ and
$$\int_{\goth X}m\otimes f(x)\text{\rm d}\mu(x)=m\otimes\int_{\goth X}
f(x)\text{\rm d}\mu(x)\text{\ in\ }N\otimes M\otimes\mu(f).$$
\noindent $(c)$ Let $f\in\Cal C_0(\goth X,M)$ and let $\phi:M\rightarrow
N$ be a homomorphism. Then $\phi\circ f$ is in $\Cal C_0(\goth X,N)$,
$\mu(\phi\circ f)\subseteq\mu(f)$ and
$$\int_{\goth X}\phi\circ f(x)\text{\rm d}\mu(x)=\phi\left(\int_{\goth X}
f(x)\text{\rm d}\mu(x)\right)\text{\ in\ }N\otimes\mu(f).$$
\noindent $(d)$ Assume that $f\in\Cal C_0(\goth X,\Bbb R)$ is locally constant and
$\mu\in\Cal M(\goth X,\Bbb R)$ is a Borel measure on $\goth X$, which is a positive measure on
the Borel sets of $\goth X$. Then the image of the integral of $f$ on $\goth X$ with respect
to $\mu$ under the homomorphism
$\Bbb R\otimes\mu(f)\rightarrow\Bbb R$ induced by the product is the usual
Lebesgue integral of $f$ on $\goth X$ with respect to $\mu$.\ $\square$
\endproclaim
\definition{Definition 5.3} Let $\Bbb R$ be a closed subfield of $\Bbb C$ and let $K\subset
\Bbb P^1(\Bbb R)$ be a non-empty compact subset. For every $\rho\in GL_2(\Bbb R)$ and
$z\in\Bbb P^1(\Bbb C)$ let $\rho(z)$ denote the image of $z$ under the M\"obius transformation
corresponding to $\rho$. Let moreover $D(\rho)$ denote the open disk
$$D(\rho)=\{z\in\Bbb P^1(\Bbb C)|1<|\rho^{-1}(z)|\}.$$
Let $\Cal D$ denote the set of open disks of the form
$D(\rho)$ where $\rho\in GL_2(\Bbb R)$. For each $D\in\Cal D$ let $D(K)$ denote 
$D\cap K$. Let $\Cal P(K)$ denote those subsets $S$ of
$\Cal D$ such that the sets $D$, $D\in S$ are pair-wise disjoint and the union of the set
$D(K)$, $D\in S$ form a partition of $K$. For each $S\in\Cal P(K)$ let $\Omega(S)$ denote
the unique connected rational subdomain defined over $\Bbb R$ with the property
$\partial\Omega(S)=S$. Let $\Omega_K$ denote the complement of $K$ in $\Bbb P^1(\Bbb C)$. Then
$\Omega_K$ is equipped naturally  with the structure of a rigid analytic space over $\Bbb R$
such that the open subsets $\Omega(S)$, $\in\Cal P(K)$ form an admissible cover by affinoid
subdomains. In particular a function $f:\Omega_K\rightarrow\Bbb R$ is holomorphic if the
restriction of $f$ onto $\Omega(S)$ is holomorphic for every $S\in\Cal P(K)$. Let
$\Cal O(\Omega_K)$ and $\Cal M(\Omega_K)$ denote the $\Bbb R$-algebra of holomorphic functions
and the field of meromorphic function of $\Omega_K$, respectively. The latter is of course the
quotient field of the former. We define $K_2(\Omega_K)$ as the intersection of the kernels of
all the tame symbols $\{\cdot,\cdot\}_x$ inside $K_2(\Cal M(\Omega_K))$ where $x$ runs through
the set $\Omega_K$. 
\enddefinition
\proclaim{Lemma 5.4} For each $k\in K_2(\Omega_K)$ there is a unique measure
$\{k\}\in\Cal M_0(K,\Bbb R^*)$ such that $\{k\}(D(K))=
\{k|_{\Omega(S)}\}_D$ for every $S\in\Cal P(K)$ and $D\in S$.
\endproclaim
\definition{Proof} First we are going to show that every open cover $\Cal U$ of $K$ has a
subordinate cover of the form $D(K)$, $D\in S$ where $S\in\Cal P(K)$. By the compactness of $K$
there is a finite cover $\Cal V$ of the form $D(K)$, $D\in I$ subordinate to $\Cal U$ where
$I\subset\Cal D$ is a finite set. We may assume that the union $\cup_{D\in I}D$ is not equal
to $\Bbb P^1(\Bbb C)$ by refining the cover $\Cal V$ further. Then any two disks in $I$ are
either disjoint or equal, hence the claim is now clear. The same argument
works for any compact and open subset $L$ of $K$. When we apply it to the one element cover
of $L$ we get that $L$ can be written as the pair-wise disjoint union of sets of the form
$D(K)$. In particular the measure $\{k\}$ is unique, if it exists. In order to prove that
$\{k\}(L)$ is well-defined we have to show that the product
$$\prod_{D\in I}\{k|_{\Omega(S)}\}_D\in\Bbb R^*$$
is independent of the choice of $I$ and $S$ for every $S\in\Cal P(K)$ and $I\subseteq S$
such that $L=\cup_{D\in I}D(K)$. Let $T\in\Cal P(K)$ and $J\subseteq T$ be another pair
such that $L=\cup_{D\in J}D(K)$. Then there is a $V\in\Cal P(K)$ and a $M\subseteq V$ such that
$\emptyset\neq\Omega(V)$ contains $\Omega(S)\cup\Omega(T)$ and $L=\cup_{D\in M}D(K)$. Clearly we
only have to show that
$$\{k|_{\Omega(U)}\}_D=\prod\Sb E(K)\subseteq D(K)\\E\in V\endSb\{k|_{\Omega(V)}\}_E$$
for every $D\in U$ where $U$ is either $S$ or $T$. By our assumptions we have either
$E(K)\subseteq D(K)$ or $E\subset\Omega_K$ for every $E\subseteq D$ with $E\in V$. As
$\{k|_{\Omega(V)}\}_E=1$ for disks of the latter type the equality above follows from the
invariance theorem. Finally note that the product $\{k\}(K)=\prod_{D\in S}
\{k|_{\Omega(S)}\}_D$ is equal to
one for every $S\in\Cal P(K)$ by $(iii)$ of Theorem 3.2 hence $\{k\}$ is indeed an
element of $\Cal M_0(K,\Bbb R^*)$.\ $\square$
\enddefinition
\definition{Definition 5.5} Now let $K$ and $L$ be two non-empty compact subsets of
$\Bbb P^1(\Bbb R)$ and assume that a non-constant holomorphic map $h:\Omega_K\rightarrow
\Omega_L$ of rigid analytic spaces over $\Bbb R$ are given. Let
$h^*:K_2(\Cal M(\Omega_L))\rightarrow K_2(\Cal M(\Omega_K))$ be the pull-back homomorphism
induced by $h$. By restriction it induces a homomorphism $K_2(\Omega_L)\rightarrow
K_2(\Omega_K)$, also denoted by $h^*$. For every abelian topological group $M$ and
compact Hausdorff topological space $\goth X$ let $\widetilde{\Cal C}_0(\goth X,M)$ 
denote the quotient of $\Cal C_0(\goth X,M)$ by the group of $M$-valued constant
functions. The integration introduced in Definition 5.1 induces a
canonical identification between Hom$(\widetilde{\Cal C}_0(\goth X,\Bbb Z),R)$ and
$\Cal M_0(\goth X,R)$ for every $R$ commutative group when $\Bbb Z$ is discrete. We are going
to define a homomorphism $h_*:\widetilde{\Cal C}_0(K,\Bbb Z)\rightarrow
\widetilde{\Cal C}_0(L,\Bbb Z)$, where $\Bbb Z$ is equipped with the discrete topology, as
follows. Given an element $\widetilde f\in\widetilde{\Cal C}_0(K,\Bbb Z)$ first choose one of
its representatives $f\in\Cal C_0(K,\Bbb Z)$. Then choose an $S\in\Cal P(K)$ such that $f$ is
equal to a constant $f(D)$ on the set $D(K)$ for every $D\in S$. Then there is a
$T\in\Cal P(L)$ such that $h(\Omega(S))\subseteq\Omega(T)$. Choose a
$\widetilde g=\sum_{E\in T}\widetilde g(E)\in\Bbb Z\partial\Omega(T)$ which represents
$$H_1(h)(\sum_{D\in S}f(D)D)\in H_1(\Omega(T)).$$
Let $g\in\Cal C_0(L,\Bbb Z)$ be the function given by the rule $g(z)=\widetilde g(E)$ for
every $z\in E(K)$ and $E\in T$. We define $h_*(\widetilde f)$ as the class of $g$ in
$\widetilde{\Cal C}_0(L,\Bbb Z)$. One may see that $h_*$ is a well-defined homomorphism
in the same way we proved the lemma above. Let $h_*:\Cal M_0(L,R)\rightarrow\Cal M_0(K,R)$
be the homomorphism induced by this $h^*$ via the duality described above. The following
proposition is an immediate consequence of the definitions and the invariance theorem:
\enddefinition
\proclaim{Proposition 5.6} We have
$$\{h^*(k)\}=h^*(\{k\})\in\Cal M_0(K,\Bbb R^*)$$
for every $k\in K_2(\Cal M(\Omega_L))$.\ $\square$
\endproclaim
\definition{Definition 5.7} For any graph $G$ let $\Cal V(G)$ and $\Cal
E(G)$ denote its set of vertices and edges, respectively. In this paper
we will only consider such oriented graphs $G$ which are equipped with an involution $\overline{\cdot}:\Cal E(G)\rightarrow\Cal E(G)$ such that for each edge $e\in\Cal E(G)$ the original and terminal vertices of the edge $\overline e\in\Cal E (G)$ are the terminal and
original vertices of $e$, respectively. The edge $\overline e$ is called the edge $e$ with reversed orientation.  Let $R$ be a
commutative group. A function $\phi:\Cal E(G)\rightarrow R$ is called a harmonic $R$-valued cochain, if it satisfies the following conditions:
\roster
\item"$(i)$" We have:
$$\phi(e)+\phi(\overline e)=0\text{\ }(\forall e\in\Cal E(G)).$$
\item"$(ii)$" If for an edge $e$ we introduce the notation $o(e)$ and
$t(e)$ for its original and terminal vertex respectively, then for all but finitely many $e\in\Cal E(G)$ with $o(e)=v$ we have $\phi(e)=0$ and
$$\sum\Sb e\in\Cal E(G)\\o(e)=v\endSb\phi(e)=0\text{\ }
(\forall v\in\Cal V(G)),$$
where by our assumption the sum above is well-defined. 
\endroster
We denote by $\Cal H(G,R)$ the group of $R$-valued harmonic
cochains on $G$.
\enddefinition
\definition{Definition 5.8} A path $\gamma$ on an oriented graph $G$ is a sequence of edges
$$\{e_1,e_2,\ldots,e_n,\ldots\}\in\Cal E(G)$$
indexed by the set $I$ where $I=\Bbb N$ or $I=\{1,\ldots,m\}$ for some
$m\in\Bbb N$ such that $t(e_i)=o(e_{i+1})$ for every $i$, $i+1\in I$. We say that
$\gamma$ is an infinite path or a finite path whether we are in the first or in the second
case, respectively. We say that a path $\{e_1,\ldots,e_n,\ldots\}$ indexed by the set
$I$ on an oriented graph $G$ is without backtracking if $\overline{e_i}\neq e_{i+1}$ for every
for every $i$, $i+1\in I$. An oriented graph $G$ is called a tree if for every pair of
different vertices $v$ and $w\in\Cal V(G)$ there is a unique finite path $\{e_1,\ldots,e_n\}$
without backtracking such that $o(e_1)=v$ and $t(e_n)=w$. Recall that a half-line $\gamma$ on
an oriented graph $G$ is an infinite path without backtracking. We say that two half-lines
on an oriented graph are equivalent if they only differ in a finite graph. We define the
set Ends$(G)$ of ends of a tree $G$ as the equivalence classes of half-lines of $G$. There is a
natural topology on Ends$(G)$ given by the sub-basis $G_e$, $e\in\Cal E(G)$ where $G_e$
consists of the equivalence classes of half-lines of the form $\{e_1,e_2,\ldots,e_n,\ldots\}$
with the property $e_1=e$. 
\enddefinition
\definition{Definition 5.9} By slightly extending the usual terminology we will call a scheme
$C$ defined over a field a curve if it is reduced, locally of finite type and of dimension one.
A curve $C$ is said to have normal crossings if every singular point of $C$ is an ordinary
double point in the usual sense. For any curve $C$ with normal crossings let $\widetilde C$
denote its normalization, and let $\widetilde S(C)$ denote the pre-image of the set $S(C)$ of
singular points of $C$. We denote by $\Gamma(C)$ the oriented graph whose set of vertices is
the set of irreducible components of $\widetilde C$, and its set of edges is the set
$\widetilde S(C)$ such that the initial vertex of any edge $x\in\widetilde S(C)$ is the
irreducible component of $\widetilde C$ which contains $x$ and the terminal vertex of $x$ is
the irreducible component which contains the unique other element $\overline x$ of
$\widetilde S(C)$ which maps with respect to the normalization map to the same singular point
as $x$. The map $x\mapsto\overline x$ is an involution
$\overline{\cdot}:\Cal E(\Gamma(C))\rightarrow\Cal E(\Gamma(C))$ of the type described in
Definition 5.7.
\enddefinition
\definition{Definition 5.10} Let $\Cal O$, $\bold f$ denote the valuation ring of
$\Bbb R$ and its residue field, respectively. Let $\goth U$ be an admissible formal scheme of
dimension one over $\Cal O$ and let $U$ denote the rigid analytic space we get from $\goth U$
by applying Raynaud's functor (for its definition see [2]). Let $\goth U_0$ denote the special
fiber of $\goth U$ and assume that the curve $\goth U_0$ over $\bold f$ is totally degenerate.
The latter means that $\goth U_0$ has normal crossings and its irreducible components are
smooth projective rational curves over $\bold f$. Assume that $U$ is biholomorphic to 
$\Omega_{\partial U}$ for some non-empty $\partial U\subseteq\Bbb P^1(\Bbb R)$. In this case
the graph $\Gamma(\goth U_0)$ is a tree and the topological space Ends$(\Gamma(\goth U_0))$
is canonically homeomorphic to $\partial U$ (see [6]). We will use this identification without
further notice. For every element $k\in K_2(U)$ let
$\{k\}$ denote the function $\{k\}:\Cal E(\Gamma(\goth U_0))\rightarrow\Bbb R^*$ which is
given by the rule $\{k\}(e)=\{k\}(\Gamma(\goth U_0)_e)$ for every edge
$e\in\Cal E(\Gamma(\goth U_0))$ where we  use the notation of Definition 5.8 and the symbol
$\{k\}$ on the right hand side of the equation above denotes the measure we introduced in
Lemma 5.4.
\enddefinition
\proclaim{Lemma 5.11} The function $\{k\}$ lies in $\Cal H(\Gamma(\goth U_0),\Bbb R^*)$.
\endproclaim
\definition{Proof} The claim is purely graph-theoretical in nature. In fact for every
tree $\Cal T$, commutative group $R$ and measure $\mu\in\Cal M_0(\text{Ends}(\Cal T),R)$ the
function $c(\mu):\Cal E(\Cal T)\rightarrow R$, given by the rule $c(\mu)(e)=\mu
(\Cal T_e)$ for every edge $e\in\Cal E(\Cal T)$, is an $R$-valued harmonic cochain. Fix a
vertex $v\in\Cal V(\Cal T)$: then every end of $\Cal T$ has a unique representative of the
form $\{e_1,e_2,\ldots,e_n,\ldots\}$ with the property $o(e_1)=v$. Now it is clear that the
sets $\Cal T_e$, $o(e)=v$ form a pair-wise disjoint partition of Ends$(\Cal T)$. Therefore
$c(\mu)$ satisfies property $(ii)$ of Definition 5.7. Similarly property $(i)$ of
Definition 5.7 follows from the fact that, given an edge $e\in\Cal E(\Cal T)$, every end of
$\Cal T$ has a unique representative of the form $\{e_1,e_2,\ldots,e_n,\ldots\}$ such that
either $e_1=e$ or $e_1=\overline e$.\ $\square$
\enddefinition
\definition{Definition 5.12} By the lemma above we have constructed a regulator
$$\{\cdot\}:K_2(U)\rightarrow\Cal H(\Gamma(\goth U_0),\Bbb R^*).$$
We are going to recall the definition of a similar invariant
$$\text{Reg}:K_2(U)\rightarrow \Cal H(\Gamma(\goth U_0),\Bbb Q)$$
which is perhaps best to call the tame regulator. Let $\nu:\Bbb R^*\rightarrow\Bbb Q$ be a
valuation induced by the absolute value $|\cdot|$. The normalization map identifies the
irreducible components of $\goth U_0$ and its normalization which we will use without further
notice.  For every vertex $v\in\Cal V(\Gamma(\goth U_0))$ let $\goth U_v$ denote open affine
subvariety of $\goth U_0$ consisting of irreducible component $v$ with all singular points
removed. Let $\goth U_v$ denote also the unique open affine formal subscheme of $\goth U$
whose fiber in $\goth U_0$ is equal to $\goth U_v$. Let $U_v$ denote the open affinoid of the
rigid analytic space $U$ we get by applying Raynaud's functor to the admissible formal scheme
$\goth U_v$. Then $U_v$ is a connected rational subdomain of $\Bbb P^1(\Bbb C)$ via the
embedding of $U$ into the latter which is isomorphic to $U(S)$ for some finite subset $S\subset
\bold f$. As we saw in Definition 4.5 there is a valuation on the
field $\Cal M(U_v)$, hence on the field  $\Cal M(U)$, whose restriction to $\Bbb R$ is $\nu$.
Let $T_v$ denote the corresponding generalized tame symbol from $K_2(U)$ into the
multiplicative group of the function field of $\goth U_v$ tensored with $\Bbb Q$. For every
$k\in K_2(U)$ and every $e\in\Cal E(\Gamma(\goth U_0))$ let Reg$(k)(e)\in\Bbb Q$ denote the
valuation of $T_{o(e)}(k)$ at the image of $e$ with respect to the normalization map. It is 
not difficult to see that Reg$(k)$ is a harmonic cochain but this fact also follows from the
following result:
\enddefinition
\proclaim{Theorem 5.13} For every $k\in K_2(U)$ we have:
$\text{\rm Reg}(k)=\nu(\{k\})$.
\endproclaim
\definition{Proof} For every $v\in\Cal V(\Gamma(\goth U_0))$ there is a bijection $b_v$
from the set
$$B_v=\{e\in\Cal E(\Gamma(\goth U_0))|o(e)=v\}$$
to the set $\partial U_v$ such
that $\nu(\{k\}_{b_v(e)})=\text{Reg}(k)(e)$ by Proposition 4.6. Since for every $e\in B_v$
we have $\Gamma(\goth U_0)=b_v(e)(\partial U)$ the claim is now obvious.\ $\square$
\enddefinition
\definition{Definition 5.14} Let $X$ be a geometrically connected regular projective curve defined over the field $\Bbb R$ and let $\Cal R(X)$ denote the field of rational functions of the curve $X$. For any $x\in X(\Bbb C)$ and any two non-zero functions $f$, $g\in\Cal R(X)$ let $\{f,g\}_x$ denote the tame symbol of the pair $(f,g)$ at $x$. We define $K_2(X)$ as the intersection of the kernels of all the tame symbols $\{\cdot,\cdot\}_x$ inside $K_2(\Cal R(X))$ where $x$ runs through the set $X(\Bbb C)$. By the usual abuse of notation let $X$ denote also the rigid analytic variety associated to the
projective curve $X$ as well.
\enddefinition
\definition{Definition 5.15} Recall that $X$ is called a Mumford curve if there is a flat, projective, regular and semistable scheme $\goth X$ over the spectrum of $\Cal O$ whose
generic fiber over $\Bbb R$ is isomorphic to $X$ and whose special fiber $\goth X_0$ over $\bold f$ is totally degenerate. Let $\goth p:\goth U\rightarrow\goth X$ be the universal cover of $\goth X$ in the category of admissible formal schemes and let $\Gamma$ denote the group of deck transformations of the cover $\goth p$. According to [6] the formal scheme $\goth U$ is of the type considered in Definition 5.10, at least under some assumptions on the base field $\Bbb R$. Let $\goth p_0:\goth U_0\rightarrow\goth X_0$ denote
the special fiber of $\goth p$ over $\bold f$ and let $p:U\rightarrow X$ denote the
map of rigid analytic spaces we get by applying Raynaud's functor to $\goth p$. The map
$\goth p_0$ induces a map of oriented graphs
$\Gamma(\goth U_0)\rightarrow\Gamma(\goth X_0)$ which in turn induces a map:
$$\goth p_0^*:\Cal H(\Gamma(\goth X_0),\Bbb R^*)\rightarrow
\Cal H(\Gamma(\goth U_0),\Bbb R^*)^{\Gamma}$$
where the superscript $\Gamma$ denotes the subgroup of $\Gamma$-invariant harmonic cochains.
The natural action of $\Gamma$ on the graph $\Gamma(\goth U_0)$ is proper and free therefore $\goth p_0^*$ is in fact an isomorphism. By the invariance theorem the regulator of the pull-back $p^*(k)$ of any element $k\in K_2(X)$ with respect to $p$
lies in $\Cal H(\Gamma(\goth U_0),\Bbb R^*)^{\Gamma}$. Hence we may define the rigid analytic
regulator for $X$ as a map:
$$\{\cdot\}:K_2(X)\rightarrow\Cal H(\Gamma(\goth X_0),\Bbb R^*)$$ 
given by the rule $\{k\}=(\goth p_0^*)^{-1}(\{p^*(k)\})$ for every $k\in K_2(X)$.
\enddefinition
\definition{Example 5.16} A case of particular interest is when $X=E$ is a Tate elliptic curve. In this case the regulator is uniquely determined by its value on any of the edges of the
reduction graph of a minimal model of the elliptic curve $E$ over Spec$(\Bbb C^0)$ so it is its really a homomorphism $\{\cdot\}:K_2(E)\rightarrow\Bbb C^*$, well-defined up to sign. It
is very easy to give an explicit description of this homomorphism in general using the Tate
uniformization of the elliptic curve. Recall that an elliptic curve $E$ defined over $\Bbb C$
is a Tate curve if its $j$ invariant $j(E)$ is not an element of $\Bbb C^0$. Under this
assumption there is a rigid analytic Tate uniformization $u:\Bbb C^*\rightarrow E$. The
pull-back $u^*(k)$ of every $k\in K_2(E)$ as an element of $K_2(\Cal R(E))$ in
$K_2(\Cal M(U))$ for any $U\subset\Bbb C^*$ connected rational subdomain lies in $K_2(U)$
hence the regulator $\{u^*(k)\}_D\in\Bbb C^*$ is well-defined for every $D\in\partial U$. The
value of this regulator $\{u^*(k)\}_D\in\Bbb C^*$ does not depend on the choice of $U$ or $D$
if the disk $D$ contains $0$ by the homotopy invariance of the regulator. This value is the
regulator $\{k\}$ of the element $k\in K_2(E)$.
\enddefinition
Next we present a purely analytical proof of Weil's Reciprocity Law for Tate elliptic curves.
\proclaim{Theorem 5.17} Let $E$ be an elliptic curve defined over $\Bbb C$ such that its
$j$-invariant $j(E)\notin\Bbb C^0$ and let $f$, $g$ be two non-zero rational functions on $E$.
Then the product of all tame symbols of the pair $(f,g)$ is equal to 1:
$$\prod_{x\in E(\Bbb C)}\{f,g\}_x=1.$$
\endproclaim
\definition{Proof} This argument can be generalized to Mumford curves using the concept of a
fundamental domain for a Schottky group, but for the sake of simplicity we present the
argument for Tate curves only. As we already noted there is a rigid analytic Tate uniformization
$u:\Bbb C^*\rightarrow E$ with Tate period $t\in\Bbb C$ such that $|t|>1$. Let $f$ and $g$ also
denote the pull-back of these functions to $\Bbb C^*$ via $u$ by abuse of notation.
Then the restriction of $f$, $g$ to the annulus $A(1,|t|)$ is
meromorphic. Let $S\subset A(1,|t|)$ be a finite set such that the
functions $f$ and $g$ don't have a zero or a pole on $A(1,|t|)-S$. For
each $s\in S$ take an open disk $D_s\subset A(1,|t|)$ containing $s$
such that these disks are pair-wise disjoint. The set $Y=A(1,|t|)-
\bigcup_{s\in S}D_s$ is a connected rational subdomain. If these disks are
sufficiently small claim $(iii)$ of Theorem 3.2 reads as follows:
$$\{f,g\}_{D(0,1)}\cdot\{f,g\}_{D(\infty,|t|^{-1})}\cdot\prod_{s\in S}
\{f,g\}_s=1.$$
But the functions $f$, $g$ are periodic with multiplicative period $t$, so the regulators $\{f,g\}_{D(0,1)}$ and $\{f,g\}_{D(\infty,|t|^{-1})}^{-1}=\{f,g\}_{D(0,|t|)}$ are equal, because they
only depend on the restrictions of $f$ and $g$ to the sets $\partial
D(0,1)$ and $\partial D(0,|t|)$, respectively. Hence the product of
the first two terms in the equation above is one, and the claim
follows.\ $\square$
\enddefinition

\heading 6. The rigid analytic regulator on Drinfeld's upper half plane
\endheading

\definition{Example 6.1} Let $\Bbb R$ denote again a closed subfield of $\Bbb C$ and assume
that the valuation on $\Bbb R$ induced by $|\cdot|$ is discrete. Also assume that the residue
field of $\Bbb R$ is a finite field $\Bbb F_q$ and let $\Cal O$ denote the valuation ring of
$\Bbb R$. Let $\Omega$ denote the rigid analytic upper half plane, or Drinfeld's
upper half plane over $\Bbb R$. It is the rigid analytic space $\Omega_K$ introduced in
Definition 5.3 in the special case when $K=\Bbb P^1(\Bbb R)$. In particular the set of points
of $\Omega$ is $\Bbb P^1(\Bbb C)-\Bbb P^1(\Bbb R)$, denoted also by $\Omega$ by abuse of
notation. We can give a very simple description of the regulator of every $k\in K_2(\Omega)$
as follows. By the invariance theorem the value  $\{k\}(\rho)=\{k|_{\Omega(S)}\}_{D(\rho)}$,
where $\rho\in GL_2(\Bbb R)$ and $D(\rho)\in S\in\Cal P$, is independent of the choice of $S$.
We define the regulator $\{k\}:GL_2(\Bbb R)\rightarrow\Bbb C^*$ of $k$ as the function given
by this rule. The assignment $k\rightarrow\{k\}$ is $GL_2(\Bbb R)$-equivariant by the
invariance theorem.
\enddefinition
\definition{Definition 6.2} We say that an additive submodule $A\subset\Bbb R$
is a lattice if it is discrete and the quotient $A\backslash\Bbb R$ is compact. Let $\Gamma(A)$ denote the following subgroup:
$$\Gamma(A)=\{\left(\matrix1&a\\0&1\endmatrix\right)\in GL_2(\Bbb R)|a\in A\}.$$
The subgroup $\Gamma(A)$ stabilizes the point $\infty$ on the
projective line via the M\"obius action. Also
note that $\Gamma(A)$ leaves the set
$$\Omega_c=\{z\in\Omega|c<|z|_i\}$$
invariant for any positive $c\in|\Bbb C|$ where $|z|_i=\inf_{x\in\Bbb R}|z-x|$
is the imaginary absolute value. Let $\mu$ be a Haar measure on the additive group of the non-archimedean field $\Bbb R$. This measure induces another measure on the quotient group $A\backslash\Bbb R$ which will be denoted my the same symbol. We may normalize $\mu$ such that $\mu(A\backslash\Bbb R)=1$. In this case $\mu$ will take only rational values. We may and we will assume that the absolute value $|\cdot|$ on $\Bbb R$ is normalized such that $\mu(y\Cal O)=|y|\mu(\Cal O)$ for every $y\in\Bbb R^*$. If $k\in K_2(\Omega)$ is a $\Gamma(A)$-invariant element then the regulator $\{k\}:GL_2(\Bbb R)\rightarrow\Bbb C^*$ is also invariant with respect to the left regular action of $\Gamma(A)$. Moreover he regulator $\{k\}$ is left invariant by multiplication on the right by a compact, open subgroup of $GL_2(\Bbb R)$ hence the integrand of the integral
$$\{k\}_{\infty}=\int_{A\backslash \Bbb R}\{k\}\left(\matrix
1&x\\0&1\endmatrix\right)d\mu(x)\in\Bbb R^*\otimes\mu(\{k\}(\left(\matrix
1&\cdot\\0&1\endmatrix\right))$$
is locally constant and the integral itself is well-defined. The modulus above is a subset of $\Bbb Q$ as we already remarked.
\enddefinition
\definition{Definition 6.3} Let $A_0$ denote the intersection $A\cap\Cal O$ and let $p$ denote the characteristic of the residue field $\Bbb F_q$. The set $A_0$ is finite and it is also a vector space over $\Bbb F_p$ so its cardinality is a power of $p$. Let $e_A(z):\Omega\rightarrow\Bbb C^*$ denote the
classical Carlitz-exponential:
$$e_A(z)=z\prod_{0\neq\lambda\in A}\left(1-{z\over\lambda}\right).$$
It is well known (see for example 2.7 of [5], page 44-45) that the function $e_A^{-1}$ is $\Gamma(A)$-invariant and it is a biholomorphic map between the quotient $\Gamma(A)\backslash\Omega_c$ and a small open disk around $0$ punctured at $0$ for a sufficiently large $c$. We say that a $\Gamma(A)$-invariant meromorphic function $u$ on $\Omega$ is
meromorphic at $\infty$ if the composition of $u$ and the inverse of the
biholomorphic map $e_A^{-1}$ is meromorphic at $0$ for some (and hence all) such $c$ number. In this case we can speak about its value, order of zero or
order of pole at $\infty$. Let $\Cal M^A(\Omega)$ denote the field of $\Gamma(A)$-invariant meromorphic functions on $\Omega$ meromorphic at $\infty$. Let $K^A_2(\Omega)$ denote the intersection $K_2(\Omega)\cap K_2(\Cal M^A(\Omega))$. For every $k\in K_2(\Cal M^A(\Omega))$ we may speak about its tame symbol at $\infty$ in the sense introduced above.
\enddefinition
\proclaim{Theorem 6.4} For each element $k\in K_2^A(\Omega)$ we have $\mu(\{k\}(\left(\smallmatrix 1&\cdot\\0&1\endsmallmatrix\right)))
\subseteq\Bbb Z$ and the integral $\{k\}_{\infty}$  is equal to the tame symbol of $k$ at $\infty$ multiplied by $|A_0|$.
\endproclaim
\definition{Proof} For every positive $c\in|\Bbb C|$ let $\Cal D_c$ denote the set of those disks $D\in\Cal D$ such that there is an element $S\in\Cal P$ such that $D\in S$ and $\Omega(S)\subset\Omega_c$. The set $\Omega_c$ has the structure of a rigid analytic space such that a function $f:\Omega_c\rightarrow\Bbb C$ is holomorphic if and only if the restriction of $f$ onto $\Omega(S)$ is holomorphic for every $S\in\Cal P$ whenever $\Omega(S)\subset\Omega_c$. Let $\Cal M(\Omega_c)$ denote the field of meromorphic functions on $\Omega_c$. For each $k\in K_2(\Cal M(\Omega))$ the value $\{k\}(D)=\{k|_{\Omega(S)}\}_{D}$, where $D\in\Cal D_c$, $D\in S\in\Cal P$ and $\Omega(S)\subset\Omega_c$, is independent of the choice of $S$ and defines a function $\{k\}:\Cal D_c\rightarrow\Bbb C^*$. For every $y\in\Bbb R^*$ and $x\in\Bbb R$ we have
$$D(\left(\matrix y&x\\0&1\endmatrix\right))=\{z\in\Bbb C||z-x|>|y|\},$$
so $D(\left(\smallmatrix y&x\\0&1\endsmallmatrix\right))\in\Cal D_c$ if $|y|>c$. In particular the integral
$$\{k\}^y_{\infty}=\int_{A\backslash \Bbb R}\{k\}(D(\left(\matrix
y&x\\0&1\endmatrix\right)))d(|y|^{-1}\mu)(x)\in\Bbb R^*
\otimes|y|^{-1}\mu(\{k\}(\left(\matrix y&\cdot\\0&1\endmatrix\right))$$
is well-defined for every $y\in\Bbb R^*$ and $k\in K_2(\Cal M(\Omega))$ when $|y|>c$ and $k$ can be represented as the linear combination of symbols of $\Gamma(A)$-invariant meromorphic functions on $\Omega_c$. Under this notation we have $\{k\}_{\infty}=\{k\}_{\infty}^1$. We say that the $\Gamma(A)$-invariant meromorphic function $u$ on $\Omega_c$ is
meromorphic at $\infty$ if it satisfies the same condition as we demanded for the elements of $\Cal M^A(\Omega)$ in Definition 6.3. Let $\Cal M^A(\Omega_c)$ denote the field of $\Gamma(A)$-invariant meromorphic functions on $\Omega_c$ meromorphic at $\infty$. Let $K^A_2(\Omega_c)$ denote the intersection $K_2(\Omega_c)\cap K_2(\Cal M^A(\Omega_c))$.
\enddefinition
\proclaim{Lemma 6.5} For each element $k\in K_2^A(\Omega_c)$ we have $|y|^{-1}\mu(\{k\}(\left(\smallmatrix y&\cdot\\0&1\endsmallmatrix\right)))
\subseteq\Bbb Z$ and $\{k\}^y_{\infty}$ does not depend on the choice of
$y\in\Bbb R^*$ where $|y|>c$.
\endproclaim
\definition{Proof} Fix a uniformizer $\pi\in\Bbb R$. As
$$D(\left(\matrix y&x\\0&1\endmatrix\right))=D(\left(\matrix y&x\\0&1\endmatrix\right)\cdot\left(\matrix u&v\\0&1\endmatrix\right))$$
for every $u\in\Cal O^*$ and $v\in\Cal O$, we may assume that $y=\pi^m$ for some integer $m$ using this identity when $v=0$. On the other hand the  identity above also implies when $v=1$ that the integrand of the integral $\{k\}^y_{\infty}$ is translation-invariant with respect to the group $y\Cal O$. Since the measure of the projection of this group into $A\backslash O$ with respect to the measure $|y|^{-1}\mu$ is an integer we get that $|y|^{-1}\mu(\{k\}(\left(\smallmatrix y&\cdot\\0&1\endsmallmatrix\right)))
\subseteq\Bbb Z$ as claimed. The function $\{k\}$ satisfies the identity:
$$\{k\}(D(g))=\sum_{\epsilon\in\Bbb F_q}\{k\}(D(g\left(\matrix
\pi&\epsilon\\0&1\endmatrix\right)))$$
for all $g\in GL_2(\Bbb R)$ because the disks $D(g\left(\smallmatrix
\pi&\epsilon\\0&1\endsmallmatrix\right))$, $\epsilon\in\Bbb F_q$ give a pair-wise disjoint partition of the disk $D(g)$. An immediate consequence of this identity is the formula:
$$\split\{k\}(D(\left(\matrix \pi^m&x\\0&1\endmatrix\right)))&
q^m\mu(x+\pi^m\Cal O)=\\&\sum_{\epsilon\in\Bbb F_q}\{k\}(D(\left(\matrix \pi^{m+1}&x+\pi^m\epsilon\\0&1\endmatrix\right)))q^{m+1}
\mu(x+\pi^m\epsilon+\pi^{m+1}\Cal O)\endsplit$$
which holds for every $x\in\Bbb R$. Hence by the translation-invariance of the measure $q^{m+1}\mu$ we have:
$$\int_{A\backslash \Bbb R}\!\!\!\{k\}(D(\left(\matrix
\pi^m&x\\0&1\endmatrix\right)))d(q^m\mu)(x)
=\int_{A\backslash \Bbb R}\!\!\!\{k\}(D(\left(\matrix
\pi^{m+1}&x\\0&1\endmatrix\right)))d(q^{m+1}\mu)(x)$$
for every sufficiently small integer $m$ as claimed.\ $\square$
\enddefinition
Let us return to the proof of theorem. Choose a presentation $k=\sum_if_i
\otimes g_i$ where $f_i$, $g_i\in\Cal M^A(\Omega)$. There is a positive $c\in\Bbb R$ such that the restriction of the functions $f_i$ and $g_i$ onto the rigid analytic space $\Omega_c$ are invertible. For every element of $K_2(\Cal M^A(\Omega_c))$ we may speak about its tame symbol at $\infty$ in the sense introduced above. By Lemma 6.5 it will be sufficient to prove that the tame symbol of $k|_{\Omega_c}$ at $\infty$ is equal to the integral $\{k\}^y_{\infty}$ for some $y\in\Bbb R^*$ with the property $|y|>c$. Therefore by bilinearity it will be sufficient to prove that the tame symbol of $f\otimes g$ at infinity is equal to $\{f\otimes g\}^y_{\infty}$ for every pair of functions $f$, $g\in\Cal O^*(\Omega_c)$ because of our assumption on $c$. In fact it will be sufficient to prove this claim in the following three cases:
\roster
\item"$(i)$" the functions $f$, $g$ are  non-zero at $\infty$,
\item"$(ii)$" the function $f$ is non-zero at $\infty$ and $g=e_A$,
\item"$(iii)$" both $f$ and $g$ are equal to $e_A$.
\endroster
In the first case we need to show that $\{f\otimes g\}^y_{\infty}=1$. We are going to show that for every positive $\epsilon$ there is an $y\in\Bbb R^*$ with the property $|y|>c$ such that $\{f\otimes g\}^y_{\infty}\in U_{\epsilon}$. This is sufficient to prove the claim in the first case by Lemma 6.5. Let $f(\infty)$ and $g(\infty)\in\Bbb R^*$ denote the value of the functions $f$ and $g$ at $\infty$, respectively. Then the values of the functions $f(z)/f(\infty)$ and $g(z)/g(\infty)$ on the rigid space $\Omega_d$ are in the set $U_{\epsilon}$ for a sufficiently large $d>0$ as the set $\Omega_d$ maps to a small neighborhood of $0$ with respect to $e_A^{-1}$. Choose an element $y\in\Bbb R^*$ such that $|y|>d$. For every $x\in\Bbb R$ let $S(x)\in\Cal P$ be a set such that $D(\left(\smallmatrix y&x\\0&1\endsmallmatrix\right))\in\Omega(S(x))$ and $\Omega(S(x))\subset\Omega_d$. By our assumptions the holomorphic functions $f/f(\infty)$, $g/g(\infty)$ are in $\Cal O_{\epsilon}(\Omega_{S(x)})$ for any $x\in \Bbb R$ hence by Theorem 2.2 we have:
$$\split\{f\otimes g\}(\left(\matrix y&x\\0&1\endmatrix\right))=&
\{{f\over f(\infty)}\otimes g\}(\left(\matrix y&x\\0&1\endmatrix\right))\cdot
\{f(\infty)\otimes{g\over g(\infty)}\}(\left(\matrix y&x\\0&1\endmatrix\right))\\
\cdot&\{f(\infty\otimes g(\infty)\}(\left(\matrix y&x\\0&1\endmatrix\right))
\in U_{\epsilon}\endsplit$$
where we also used that the third factor on the right hand side is equal to $1$. Hence the integral $\{k\}^y_{\infty}$ is an element of $U_{\epsilon}$, too, since the corresponding modulus is a subset of $\Bbb Z$.

In the second case we need to show that $\{f\otimes e_A\}_{\infty}^1
=f(\infty)^{-|A_0|}$ where $f(\infty)$ denotes again the value of the function $f$ at $\infty$. It is clear that $\{(f/f(\infty))\otimes e_A\}_{\infty}^y=1$ by repeating the argument used in the proof of the claim in the first case, therefore we may assume that $f=f(\infty)$ is constant. By the definition of the degree homomorphism we have
$$\{f(\infty)\otimes e_A\}(D)=f(\infty)^{\deg(e_A|_{\Omega(S)})(D)},$$
where $D\in\Cal D_c$, $D\in S\in\Cal P$ and $\Omega(S)\subset\Omega_c$. In fact for any $u\in\Cal O^*(\Omega_c)$ the expression $\deg(u|_{\Omega(S)})(D)$ is independent of the choice of $S$ and defines a function $\deg(u):\Cal D_c\rightarrow\Bbb Z$. (It is not difficult see that this is the van der Put logarithmic differential when the domain of definition of $u$ is $\Omega$.) Hence it will sufficient to prove that
$$\int_{A\backslash\Bbb R}\deg(e_A)d\mu(x)=-|A_0|\in\Bbb Z\otimes\Bbb Z=\Bbb Z.$$
By $(ii)$ of Theorem 3.6 we have:
$$\deg(e_A)(D(g))=-|\{\lambda\in A|\lambda\notin D(g)\}|$$
for every $g\in GL_2(\Bbb R)$ such that $\infty\in D(g)$. As
$$\infty\in D(\left(\matrix1&x\\0&1\endmatrix\right))=
\{z\in\Bbb P^1(\Bbb C)|1<|z-x|\}$$
for any $x\in \Bbb R$, we get:
$$\deg(e_A)(\left(\matrix1&x\\0&1\endmatrix\right))=-|\{\lambda\in
A||\lambda-x|\leq1\}|=
\cases-|A_0|,&\text{if $x\in A+\Cal O$,}\\
\quad\!\!0,&\text{otherwise,}\endcases$$
so the claim is now clear in the second case. In the last case we need to show that $\{e_A\otimes e_A\}_{\infty}^1=(-1)^{-|A_0|}$. Note that
$$\{f\otimes f\}_D=(-1)^{\deg(f)(D)}$$
for every $U\subset\Bbb P^1$ rational subdomain, $D\in\partial U$ boundary component and $f\in\Cal O^*(U)$ function. This is obviously true for rational functions and the general case follows from this one by approximation. Therefore we only have to show that
$$\int_{A\backslash\Bbb R}(-1)^{\deg(e_A)}d\mu(x)=
(-1)^{-|A_0|}\in\Bbb R^*\otimes\Bbb Z=\Bbb R^*.$$
It is clear from the calculations above that the integrand on the left hand side above is equal to $(-1)^{-|A|_0}$, if $x\in A+\Cal O$, and it is equal to $0$, otherwise, hence the required identity obviously holds.\ $\square$

\enddocument